\input ./style/arxiv-general.cfg
\documentclass[MSNbibl,number,citesort,secthm,seceqn,dvips]{arxbj}
\makeatletter
   \@ifpackageloaded{graphicx}{}{\usepackage{graphicx}}
\makeatother
\usepackage{upgreek}

\volume{22}
\issue{3}
\pubyear{2016}
\firstpage{1572}
\lastpage{1597}
\doi{10.3150/15-BEJ704}
\docsubty{FLA}

\makeatletter
\newtheorem{hypothesis}{Hypothesis}
\newcommand{\rrvert}{\vert}
\newcommand{\rrVert}{\Vert}
\newcommand{\llvert}{\vert}
\newcommand{\llVert}{\Vert}
\newtheorem{theorem}{Theorem}[section]
\newtheorem{proposition}[theorem]{Proposition}
\newtheorem{lemma}[theorem]{Lemma}
\newremark{remark}[theorem]{Remark}
\makeatother

\begin{document}
\begin{frontmatter}

\title{Approximation of a stochastic wave equation in dimension
three, with application to a~support theorem in H\"older
norm: The~non-stationary case}
\runtitle{Approximation and support theorem}

\begin{aug}
\author[A]{\inits{F.J.}\fnms{Francisco J.}~\snm{Delgado-Vences}\thanksref{A}\ead[label=e1]{delgadovences@gmail.com}}
\and
\author[B]{\inits{M.}\fnms{Marta}~\snm{Sanz-Sol\'e}\thanksref{B}\corref{}\ead[label=e2]{marta.sanz@ub.edu}}
\address[A]{Istituto Nazionale di Geofisica e Vulcanologia, Pisa,
Italy. \printead{e1}}
\address[B]{Facultat de Matem\`atiques, Universitat de Barcelona,
Gran Via de les Corts Catalanes, 585, E-08007 Barcelona, Spain. \printead{e2}}
\end{aug}

%
\received{\smonth{4} \syear{2014}}
%
\revised{\smonth{1} \syear{2015}}

%
\begin{abstract}
This paper is a continuation of
(\textit{Bernoulli} \textbf{20} (2014) 2169--2216) where we
prove a characterization of the support in H\"older norm of
the law of the solution to a stochastic wave equation with
three-dimensional space variable and null initial conditions. Here, we
allow for non-null initial conditions and, therefore, the solution does
not possess a stationary property in space. As in
(\textit{Bernoulli} \textbf{20} (2014) 2169--2216), the support theorem is a consequence of
an approximation result, in the convergence of probability,
of a sequence of evolution equations driven by a family of
regularizations of the driving noise. However, the method of the proof
differs from (\textit{Bernoulli} \textbf{20} (2014) 2169--2216) since arguments based on the
stationarity property of the solution cannot be used.
\end{abstract}

%
\begin{keyword}
\kwd{approximating schemes}
\kwd{stochastic wave equation}
\kwd{support theorem}
\end{keyword}
\end{frontmatter}

\section{Introduction}
\label{s1}

This paper is a continuation of \cite{Delgado--Sanz-Sole012}, where we
prove a characterization of the
topological support in H\"older norm for the law of the solution of a
stochastic wave equation with
vanishing initial conditions. Consider the stochastic partial
differential\vspace*{-2pt} equation (SPDE)
%
\begin{eqnarray}
\biggl(\frac{\partial^2}{\partial t^2} - \Delta \biggr) u(t,x) &=& \sigma \bigl(u(t,x) \bigr)
\dot M(t,x) + b \bigl(u(t,x) \bigr),
\nonumber
\\[-9pt]
\label{s1.1}
\\[-9pt]
u(0,x) &=& v_0(x), \qquad \frac{\partial}{\partial t} u(0,x) = \tilde
v_0(x),
\nonumber
\end{eqnarray}
%
where $\Delta$ denotes the Laplacian on $\mathbb{R}^3$, $T>0$ is
fixed, $t\in
(0,T]$ and $x\in\mathbb{R}^3$. The
nonlinear terms and the initial conditions are defined by functions
$\sigma, b\dvtx  \mathbb{R}\rightarrow\mathbb{R}$
and $v_0, \tilde v_0 \dvtx \mathbb{R}^3 \rightarrow\mathbb{R}$,
respectively. The
notation $\dot M(t,x)$ refers to the
formal derivative of a Gaussian random field $M$ white in the time
variable and with a correlation
in the space variable given by a Riesz kernel.
More specifically,\vspace*{-1pt}
%
\begin{equation}
\mathbb{E} \bigl(\dot M(t,x) \dot M(s,y) \bigr) = \delta_{0}(t-s)
\llvert x-y\rrvert ^{-\beta},
\end{equation}
where $\delta_{0}$ denotes\vadjust{\goodbreak} the delta Dirac measure and $\beta\in(0,2)$.

We consider a random field solution to the SPDE (\ref{s1.1}), which
means a real-valued adapted
(with respect to the natural filtration generated by the Gaussian
process $M$) stochastic process
$\{u(t,x), (t,x)\in(0,T]\times\mathbb{R}^3\}$
satisfying
%
\begin{eqnarray}
u(t,x) & =&  X^0(t,x)+\int_0^t \!
\int_{\mathbb{R}^3} G(t-s,x-y) \sigma \bigl(u(s,y) \bigr) M(\mathrm{d}s,\mathrm{d}y)
\nonumber
\\[-8pt]
\label{s1.6}
\\[-8pt]
\nonumber
&&{}+ \int_0^t \bigl[G(t-s,\cdot)\star b
\bigl(u(s,\cdot) \bigr) \bigr](x) \,\mathrm{d}s.
\end{eqnarray}
Here,
%
\begin{equation}
\label{i.c}
X^0(t,x)= \bigl[G(t)\star\tilde v_0
\bigr](x)+ \biggl[\frac{\mathrm{d}}{\mathrm{d} t} G(t)\star v_0 \biggr](x),
\end{equation}
$G(t)$ is the fundamental solution to the wave equation in dimension three,
$G(t,\mathrm{d}x)=\frac{1}{4\uppi t}\sigma_t(\mathrm{d}x)$,
where $\sigma_t(x)$ denotes the uniform surface measure on the sphere
of radius $t$ with total mass
$4\uppi t^2$ (see, e.g., \cite{Folland76}), and
the symbol ``$\star$'' denotes the convolution in the spatial argument.

The stochastic integral (also termed stochastic convolution) in (\ref
{s1.6}) is defined as a
stochastic integral with respect to a sequence of independent standard
Brownian motions
$\{W_j(s)\}_{j\in\mathbb{N}}$, as follows. Let $\mathcal{H}$ be the Hilbert
space defined by the completion
of $\mathcal{S}(\mathbb{R}^3)$, the space of rapidly decreasing
functions on
$\mathbb{R}^3$, endowed with the semi-inner
product
%
\[
\langle\varphi,\psi\rangle_\mathcal{H}= \int_{\mathbb{R}^3}\mu
(\mathrm{d}\xi)\mathcal{F} \varphi(\xi) \overline{\mathcal{F}\psi(\xi)},
\]
where $\mathcal{F}$ denotes the Fourier transform operator and $\mu
(\mathrm{d}\xi)=
\mathcal{F}^{-1}(|\xi|^{-\beta}\,\mathrm{d}\xi) =|\xi|^{\beta-3} \,\mathrm{d}\xi$. Then
%
\begin{eqnarray}
&&  \int_0^t\! \int
_{\mathbb{R}^3} G(t-s,x-y) \sigma \bigl(u(s,y) \bigr) M(\mathrm{d}s,\mathrm{d}y)
\nonumber
\\[-8pt]
\label{s1.7}
\\[-8pt]
\nonumber
&&\quad :=\sum_{j\in\mathbb{N}}\int_0^t
\bigl\langle G(t-s,x-\ast) \sigma \bigl(u(s,\ast) \bigr), e_j \bigr
\rangle_{\mathcal{H}} W_j(\mathrm{d}s),
\end{eqnarray}
where $(e_j)_{j\in\mathbb{N}}\subset\mathcal{S}(\mathbb{R}^3)$ is a
complete orthonormal basis of $\mathcal{H}$.

Assume that $\varphi\in\mathcal{H}$ is a signed measure with finite total
variation. Then, by applying \cite{kx}, Theorem~5.2
(see also \cite{Mattila95}, Lemma~12.12,  for the case of probability
measures with compact support) and a
polarization argument on the positive and negative parts of $\varphi$,
we obtain
%
\begin{equation}
\label{mattila}
\Vert\varphi\Vert_\mathcal{H}^2=C\int
_{\mathbb{R}^3}\!\int_{\mathbb{R}^3} \varphi(\mathrm{d}x) \varphi(\mathrm{d}y)
|x-y|^{-\beta}.
\end{equation}
%

For $t_0\in[0,T]$, $K\subset\mathbb{R}^3$ compact and $\rho\in
(0,1)$, we
denote by $\mathcal{C}^\rho([t_0,T]\times K)$
the space of real functions $g$ such that $\Vert g\Vert_{\rho
,t_0,K}<\infty$, where
%
\[
\Vert g\Vert_{\rho,t_0,K}:= \sup_{(t,x)\in[t_0,T]\times K} \bigl\llvert g(t,x)
\bigr\rrvert + \mathop{\sup_{(t,x),(\bar{t},\bar{x})\in[t_0,T]\times
K}}_{(t,x)\ne(\bar t,\bar x)} \frac{\llvert g(t,x)-g(\bar{t},\bar{x})\rrvert  }{(\llvert t-\bar{t}\rrvert +\llvert x-\bar{x}\rrvert )^\rho
}.
\]

Let $0<\rho^\prime<\rho$ and $\mathcal{E}^{\rho^\prime
}([t_0,T]\times K)$ be the space of H\"older
continuous functions $g$ of degree $\rho^\prime$ such that
%
\begin{equation}
O_g(\delta):=\sup_{\llvert t-s\rrvert +\llvert x-y\rrvert <\delta}\frac{\vert
g(t,x)-g(s,y)\vert}{(\llvert t-s\rrvert +\llvert x-y\rrvert )^{\rho^\prime}}
\rightarrow0, \qquad\mbox{if } \delta\to0.
\end{equation}
The space $\mathcal{E}^{\rho^\prime}([t_0,T]\times K)$ endowed with
the norm $\Vert\cdot\Vert_{\rho^\prime,t_0,K}$ is a Polish
space and the embedding $\mathcal{C}^\rho([t_0,T]\times K)\subset
\mathcal{E}^{\rho^\prime}([t_0,T]\times K)$ is compact.

Assume that the functions $\sigma$ and $b$ are Lipschitz continuous
and the initial conditions $v_0$,
$\tilde v_0$ satisfy the assumption (h2) of Theorem~\ref{sth} below.
Theorem~4.11 in \cite{Dalang--Sanz-Sole09}
along with \cite{Dalang-Quer011}, Proposition~2.6, give the existence
of a random field solution to (\ref{s1.6})
with sample paths in $\mathcal{C}^\rho([0,T]\times K)$, with $\rho
\in (0,\gamma_1\wedge\gamma_2\wedge\frac{2-\beta}{2} )$.

For any $t\in(0,T]$, let $\mathcal{H}_t=L^2([0,t]; \mathcal{H})$.
Fix $h\in\mathcal{H}_t$
and consider the deterministic evolution equation
%
\begin{eqnarray}
\Phi^h(t,x) &=& X^0(t,x) + \bigl\langle
G(t-\cdot,x-\ast)\sigma \bigl(\Phi^h(\cdot,\ast) \bigr) ,h \bigr
\rangle_{\mathcal{H}_t}
\nonumber
\\[-8pt]
\label{sm.h}
\\[-8pt]
\nonumber
&&{}+\int_0^t \mathrm{d}s \bigl[G(t-s,\cdot)\star \bigl(
\Phi^h(s,\cdot) \bigr) \bigr](x).
\end{eqnarray}
The main objective in \cite{Delgado--Sanz-Sole012} is to prove that,
in the particular case $v_0=\tilde v_0=0$,
the topological support of the law of the solution to (\ref{s1.6}) in
the space $\mathcal{E}^\rho([t_0,T]\times K)$ with
$\rho\in (0,\frac{2-\beta}{2} )$ is the closure in the
H\"older norm of the set $\{\Phi^h, h\in\mathcal{H}_T\}$, for any
$t_0>0$ (see \cite{Delgado--Sanz-Sole012}, Theorem~3.1).

The aim of this paper is to prove a partial extension of this result
allowing non-null initial conditions $v_0$, $\tilde v_0$,
but restricting to affine coefficients $\sigma$. In particular, this
will apply to the \textit{hyperbolic Anderson model}
($\sigma(x)=\lambda x$, $\lambda\ne0$). The theorem is as follows.

\begin{theorem}
\label{sth}
Assume that
\begin{enumerate}[(h2)]
\item[(h1)] the function $\sigma$ is affine and $b$ is Lipschitz continuous;
\item[(h2)] $v_0, \tilde v_0 \dvtx \mathbb{R}^3 \rightarrow\mathbb{R}$
are bounded,
$v_0\in\mathcal{C}^2(\mathbb{R}^3)$, $\nabla v_0$ is bounded,
$\tilde v_0\in\mathcal{C}(\mathbb{R}^3)$, $\Delta v_0$ and $\tilde
v_0$ are
H\"older continuous functions of degree $\gamma_1, \gamma_2$,
respectively.
\end{enumerate}
Fix $t_0>0$ and a compact set $K\subset\mathbb{R}^3$. Then the topological\vspace*{1pt}
support of the law of the solution to (\ref{s1.6})
in the space $\mathcal{E}^\rho([t_0,T]\times K)$ with $\rho\in
(0, \gamma_1\wedge\gamma_2\wedge (\frac{2-\beta}{2}
) )$
is the closure in the H\"older norm $\Vert\cdot\Vert_{\rho,t_0,K}$
of the set $\{\Phi^h, h\in\mathcal{H}_T\}$, where
$\Phi^h$ is given in (\ref{sm.h}).
\end{theorem}

After the seminal paper \cite{stroock},
an extensive literature on support theorems for stochastic differential
equations appeared (see, e.g.,
\cite{aida-kusuoka-stroock93,Gyongy--Nualart--Sanz-Sole95,Millet--Sanz-Sole94b}, and references
herein). The analysis of the uniqueness of invariant measures is one of
the motivations for the characterization
of the support of stochastic evolution equations (see \cite{Delgado--Sanz-Sole012}, Section~1, for more details).

As in \cite{Delgado--Sanz-Sole012}, Theorem~\ref{sth} will be a
corollary of a general result on approximations of
equation~(\ref{s1.6}) by a sequence of SPDEs obtained by smoothing the
noise $M$. The precise statement, given in
Theorem~\ref{ts2.1}, provides a Wong--Zakai-type theorem in H\"older
norm. It is of interest by its own. The method
relies on \cite{aida-kusuoka-stroock93}, further developed and used in
\cite{Bally--Millet-Sanz-Sole95,Gyongy--Nualart--Sanz-Sole95,Millet--Sanz-Sole94a,Millet--Sanz-Sole94b,Millet--Sanz-Sole00}. We refer the
reader to \cite{Delgado--Sanz-Sole012}, Section~1,
for a detailed description of the method for the proof of support
theorems based on approximations.

In contrast with the situation considered in \cite{Delgado--Sanz-Sole012}, the solution to (\ref{s1.6}) with non-null
initial conditions does not possess the spatial stationary property
termed \textit{$\mathcal{S}$ property} in \cite{Dalang99}.
This property is crucial in the proof of the analogue of Theorem~\ref
{ts2.2} and more precisely, in establishing the
upper bound of $L^p$ norms of increments in space when the initial
conditions are null.
The new approach to the proof of a similar upper bound when the initial
conditions do not vanish uses fractional
Sobolev norms and the classical Sobolev's embeddings (see
Proposition~\ref{pss2.1.2}). To some extend, some of the
results of this paper are a refinement and an extension of results of
\cite{Dalang--Sanz-Sole09}. Compare, for example,
Lemma~\ref{lss2.1.3} with \cite{Dalang--Sanz-Sole09}, Proposition~3.5, and Proposition~\ref{pss2.1.2} with
\cite{Dalang--Sanz-Sole09}, Theorem~4.6. Others, like Proposition~\ref
{pss2.1.7}, are crucial results to establish
the approximations. The proof of Proposition~\ref{pss2.1.7} requires
the validity of the inequality
$\Vert B(f) - B(g)\Vert_{\gamma,p,\mathcal{O}} \le C \Vert f-g\Vert
_{\gamma,p,\mathcal{O}}$ for
$B \dvtx \mathbb{R}\to\mathbb{R}$ and functions $f,g$ belonging to the
fractional Sobolev space
$W^{\gamma,p}(\mathcal{O})$ (see (\ref{frac}), (\ref{sobolevfrac})
for the definition of these spaces).
This holds when $B$ is affine and not only Lipschitz, which explains
the hypothesis on $\sigma$ in Theorem~\ref{sth}. The use of
fractional norms seems to be at the origin of
this restriction, as was noticed for
example in \cite{V--S-S09}.

The paper is organized in the following way. In Section~\ref{s2}, we
prove Theorem~\ref{ts2.1} -- a general
result on approximations of SPDEs in the convergence of probability and
in the H\"older norm -- which in
turn follow from Theorems \ref{ts2.2} and \ref{ts2.3}. As a
particular case, the characterization of the
support stated in Threorem \ref{sth} is established.
Section~\ref{sa} gathers some technical results used in the proofs.

\section{Approximations of the wave equation}
\label{s2}

As in the companion paper \cite{Delgado--Sanz-Sole012}, we consider
smooth approximations of $W$ defined as follows.
For any $n\in\mathbb{N}$, we define the partition of $[0,T]$
consisting of
the points $\frac{iT}{2^n}$, $i=0,1,\ldots,2^n$.
Denote by $\Delta_i$ the interval $ [\frac{iT}{2^n} , \frac
{(i+1)T}{2^n} )$ and by $|\Delta_i|$ its length.
We write $W_j(\Delta_i)$ for the increment $W_j(\frac{(i+1)T}{2^n})-
W_j(\frac{iT}{2^n})$, $i=0,\ldots,2^n-1$, $j\in\mathbb{N}$.
Then we define $W^n$ as the sequence whose terms are
\[
W_j^n=\int_0^\cdot
\dot{W}_j^n(s) \,\mathrm{d}s,\qquad j\in\mathbb{N},
\]
where for $j>n$, $\dot{W}_j^n=0$, and for $1\le j\le n$,
\[
\dot{W}_j^n(t)= %
\cases{
\displaystyle\sum_{i=0}^{2^n-2}
2^{n}T^{-1} W_j(\Delta_i)1_{\Delta_{i+1}}(t),
&\quad\mbox{if }$t\in \bigl[2^{-n}T,T \bigr]$,\vspace*{3pt}
\cr
0, &
\quad\mbox{if }$t \in \bigl[0,2^{-n}T \bigr[\,$. }
\]

Set
%
\begin{equation}
\label{s2.1}
w^n(t,x)=\sum_{j\in\mathbb{N}}
\dot{W}_j^n(t) e_j(x).
\end{equation}
It can be easily checked that, for any $p\in[2,\infty)$,
%
\begin{equation}
\label{s2.2}
\bigl\|w^n\bigr\|_{L^p(\Omega,\mathcal{H}_T)}\le C n^{{1}/{2}}2^{{n}/{2}}.
\end{equation}
Hence, $w^n$ belongs to $\mathcal{H}_T$ a.s.

We consider the integral equations
\begin{eqnarray}
X(t,x)&=&X^0(t,x)+\int_0^t\! \int
_{\mathbb{R}^3} G(t-s,x-y) (A+B) \bigl(X(s,y) \bigr) M(\mathrm{d}s,\mathrm{d}y)
\nonumber
\\[-8pt]
\label{s2.3}
\\[-8pt]
\nonumber
&&{}+ \bigl\langle G(t-\cdot,x-\ast)D \bigl(X(\cdot,\ast) \bigr),h \bigr
\rangle_{\mathcal{H}_t} +\int_0^t
\bigl[G(t-s,\cdot)\star b \bigl(X(s,\cdot) \bigr) \bigr](x)
\,\mathrm{d}s,
\\
X_n(t,x)&=&X^0(t,x)+\int_0^t\!
\int_{\mathbb{R}^3} G(t-s,x-y) A \bigl(X_n(s,y) \bigr)
M(\mathrm{d}s,\mathrm{d}y)
\nonumber
\\
\label{s2.4}
&&{}+ \bigl\langle G(t-\cdot,x-\ast)B \bigl(X_n(\cdot,\ast)
\bigr),w^n \bigr\rangle_{\mathcal{H}_t} + \bigl\langle G(t-\cdot,x-
\ast)D \bigl(X_n(\cdot,\ast) \bigr),h \bigr\rangle_{\mathcal{H}_t}\qquad
\\
\nonumber
&&{}+\int_0^t \bigl[G(t-s,\cdot)\star b
\bigl(X_n(s,\cdot) \bigr) \bigr](x) \,\mathrm{d}s,
\end{eqnarray}
where $h\in\mathcal{H}_T$, $w^n$ defined as in (\ref{s2.1}), $A, B,
D, b \dvtx \mathbb{R}\to\mathbb{R}$, and $X^0 \dvtx [0,T]\times\mathbb{R}^3\to
\mathbb{R}$
is the deterministic function defined in (\ref{i.c}).

For each $t\in[0,T]$, let $t_n= \max\{\underline{t}_n- 2^{-n}T, 0 \}$, with
%
\begin{equation}
\label{s2.7} \underline{t}_n = \max \bigl\{k2^{-n}T,
k=0,\ldots,2^n-1 \dvt k2^{-n}T\le t \bigr\}.
\end{equation}
By means of the following expressions, we define stochastic processes
close to $X(t_n,x)$ and $X_n(t_n,x)$,
$(t,x)\in[0,T]\times\mathbb{R}^3$, respectively:
%
\begin{eqnarray}
 X(t,t_n,x)&=&X^0(t,x)+\int
_0^{t_n} \!\!\int_{\mathbb{R}^3} G(t-s,x-y)
(A+B) \bigl(X(s,y) \bigr) M(\mathrm{d}s,\mathrm{d}y)
\nonumber
\\
\label{s2.6}
&&{}+ \bigl\langle G(t-\cdot,x-\ast)D \bigl(X(\cdot,\ast) \bigr)1_{[0,t_n]}(
\cdot ),h \bigr\rangle_{\mathcal{H}_t}
\\
&&{}+\int_0^{t_n} \bigl[G(t-s,\cdot)\star b
\bigl(X(s,\cdot) \bigr) \bigr](x) \,\mathrm{d}s,\nonumber
\\
X_n^-(t,x) &=& X^0(t,x)+\int
_0^{t_n} \!\!\int_{\mathbb{R}^3} G(t-s,x-y)
A \bigl(X_n(s,y) \bigr) M(\mathrm{d}s,\mathrm{d}y)
\nonumber
\\
&&{}+ \bigl\langle G(t-\cdot,x-\ast)B \bigl(X_n(\cdot,\ast)
\bigr)1_{[0,t_n]}(\cdot ),w^n \bigr\rangle_{\mathcal{H}_t}
\nonumber
\\[-8pt]
\label{s2.5}
\\[-8pt]
\nonumber
&&{}+ \bigl\langle G(t-\cdot,x-\ast)D \bigl(X_n(\cdot,\ast)
\bigr)1_{[0,t_n]}(\cdot ),h \bigr\rangle_{\mathcal{H}_t}
\nonumber
\\
&&{}+\int_0^{t_n} \bigl[G(t-s,\cdot)\star b
\bigl(X_n(s,\cdot) \bigr) \bigr](x) \,\mathrm{d}s.\nonumber
\end{eqnarray}

We will consider the following set of assumptions.

\renewcommand{\thehypothesis}{(H)}
\begin{hypothesis}\label{H}
\begin{enumerate}[(H2)]
\item[(H1)]  The functions $A,B,D,b \dvtx  \mathbb{R}\mapsto\mathbb{R}$ are globally
Lipschitz continuous.
\item[(H2)] $v_0, \tilde v_0 \dvtx  \mathbb{R}^3 \rightarrow\mathbb{R}$
are bounded,
$v_0\in\mathcal{C}^2(\mathbb{R}^3)$, $\nabla v_0$ is
bounded, $\tilde v_0\in\mathcal{C}(\mathbb{R}^3)$, $\Delta v_0$ and
$\tilde
v_0$ are H\"older continuous functions of degree
$\gamma_1, \gamma_2$, respectively.
\end{enumerate}

Let $\mathcal{O}$ be a bounded or unbounded open subset of $\mathbb{R}^3$,
$p\in[1,\infty)$, $\gamma\in(0,1)$. We define
%
\begin{equation}
\label{frac} \Vert g \Vert_{\gamma,p,\mathcal{O}}= \biggl(\int_{\mathcal{O}}
\mathrm{d}x \int_{\mathcal{O}} \mathrm{d}y \frac{\vert g(x)-g(y)\vert^p}{
\vert x-y\vert^{3+\gamma p}} \biggr)^{{1}/{p}}.
\end{equation}
Then we denote by $W^{\gamma,p}({\mathcal{O}})$ the Banach space
consisting of functions $g \dvtx  \mathbb{R}^3 \to\mathbb{R}$ such that
%
\begin{equation}
\label{sobolevfrac} \Vert g\Vert_{W^{\gamma,p}({\mathcal{O}})}:= \bigl(\Vert g\Vert
_{L^p(\mathcal{O})}^p + \Vert g \Vert_{\gamma,p,\mathcal{O}}^p
\bigr)^{{1}/{p}} < \infty.
\end{equation}
From \cite{Dalang--Sanz-Sole09}, Lemmas~4.2 and~4.4 and
\cite{Dalang-Quer011}, Lemma~4.2, it follows that
\textup{(H2)} implies the following.
\begin{enumerate}[(H2.2)]
\item[(H2.1)] For any $t\in[0,T]$ and any bounded domain $\mathcal
{O}\subset\mathbb{R}^3$, for any $p\in[2,\infty)$ such
that $\frac{2-\beta}{2}>\frac{3}{p}$, and for any $\gamma\in
(0,\gamma_1\wedge\gamma_2\wedge(\frac{2-\beta}{2}-
\frac{3}{p}) )$,
\[
\bigl\Vert X^0(t)\bigr\Vert_{W^{\gamma,p}(\mathcal{O})}<\infty.
\]
\item[(H2.2)] $(t,x)\mapsto X^0(t,x)$ is continuous and $\sup_{(t,x)\in[0,T]\times\mathbb{R}^3}|X^0(t,x)| <\infty$.
\end{enumerate}
\end{hypothesis}

The existence and uniqueness of a random field solution to the
equations (\ref{s2.3}), (\ref{s2.4}) is established as
in \cite{Delgado--Sanz-Sole012}, Theorem~5.1. It is proved using the
convergence of a Picard iteration scheme. For~(\ref{s2.3}), the Picard approximations converge in $L^p(\Omega)$,
uniformly in $(t,x)\in[0,T]\times\mathbb{R}^3$. For~(\ref{s2.4}), the convergence of the Picard approximations holds in
probability. It is obtained using a localization
in $\Omega$. Notice that equation (\ref{s2.4}) is more general than
(\ref{s2.3}).

The aim of this section is to prove the following theorem, which is the
analogue of \cite{Delgado--Sanz-Sole012}, Theorem~2.2,
in the context of this article.

\begin{theorem}
\label{ts2.1}
We assume Hypothesis~\textup{\ref{H}} and in addition that the function $B$ is
affine. Fix $t_0>0$ and a compact set $K\subset\mathbb{R}^3$.
Then for any $\rho\in (0,\gamma_1\wedge\gamma_2\wedge\frac
{2-\beta}{2} )$ and $\lambda>0$,
%
\begin{equation}
\label{s2.8} \lim_{n\to\infty} \mathbb{P} \bigl(\|X_n-X
\|_{\rho,t_0,K}> \lambda \bigr)=0.
\end{equation}
\end{theorem}

With a particular choice of the functions $A$, $B$ and $D$ in equations
(\ref{s2.3}), (\ref{s2.4}), this theorem yields
the characterization of the support stated in Theorem~\ref{sth} (see
the proof of Theorem~3.1 in~\cite{Delgado--Sanz-Sole012}).

The proof of Theorem~\ref{ts2.1} entails several steps. As in the
stationary case considered in \cite{Delgado--Sanz-Sole012}, the
main ingredients are \textit{local} $L^p$ estimates of increments of $X_n$
and $X$, in time and in space, and a \textit{local} $L^p$
convergence of the sequence $X_n(t,x)$ to $X(t,x)$. Here, in contrast
with \cite{Delgado--Sanz-Sole012}, \textit{local} $L^p$
estimates of increments of $X_n$ and $X$ in space are obtained via
Sobolev's embeddings.


We remind the \textit{localization} procedure introduced in \cite
{Millet--Sanz-Sole00} and also used in \cite{Delgado--Sanz-Sole012}.
For any integer $n\ge1$ and $t\in[0,T]$, define
%
\begin{equation}
\label{s2.9}
L_{n}(t)= \Bigl\{\sup_{1\le j\le n}
\sup_{0\le i\le
[2^nt T^{-1}-1]^+} \bigl\llvert W_j(
\Delta_i) \bigr\rrvert \le \alpha n^{{1}/{2}}2^{-{n}/{2}}
\Bigr\},
\end{equation}
with $\alpha>(2\ln2)^{{1}/{2}}$.
The mapping $t\mapsto L_n(t)$ is decreasing and $\lim_{n\to\infty}
\mathbb{P}(L_n(t)^c)=0$
(see \cite{Millet--Sanz-Sole00}, Lemma~2.1). It is easy to check that
%
\begin{equation}
\label{s2.10} \bigl\Vert w^n(t,\ast)1_{L_n(t)}
\bigr\Vert_{\mathcal{H}}\le C n^{{3}/{2}} 2^{{n}/{2}},
\end{equation}
and also 
\[
 \bigl\llVert w^{n}1_{L_n(t')}1_{[t,t']}
\bigr\rrVert _{\mathcal{H}_T}\le Cn^{{3}/{2}}2^{{n}/{2}} \bigl\llvert
t'-t \bigr\rrvert ^{{1}/{2}},\qquad 0\le t\le t'\le T.
\]
In particular, if $[t,t']\subset\Delta_i$ for some $i=0,\ldots
,2^n-1$, then 
\begin{equation}
\label{s2.11}
\bigl\llVert w^{n}1_{L_n(t')}1_{[t,t']}
\bigr\rrVert _{\mathcal{H}_T}\le Cn^{{3}/{2}}.
\end{equation}
%



As has been said in the \hyperref[s1]{Introduction}, the proof of Theorem~\ref{ts2.1}
will follow from Theorems \ref{ts2.2} and
\ref{ts2.3} below. These are the analogues of \cite{Delgado--Sanz-Sole012}, Theorems~2.3 and~2.4, in the context
of this article. We denote by $\Vert\cdot\Vert_p$ the $L^p(\Omega)$
norm, and for any compact set $K\subset\mathbb{R}^3$, we define
%
\begin{equation}
\label{loc} K(t)= \bigl\{x\in\mathbb{R}^3 \dvt d(x,K)\le T-t
\bigr\}, \qquad t\in[0,T],
\end{equation}
where $d$ denotes the Euclidean distance. Notice that $t\mapsto K(t)$
is a decreasing mapping.

\begin{theorem}
\label{ts2.2}
We assume Hypothesis~\textup{\ref{H}} and also that the function $B$ is affine.\vspace*{1pt} Fix
$t_0>0$ and $t_0\le t\le\bar{t}\le T$,
$x,\bar{x}\in\mathbb{R}^3$. Then, for any $p\in[1,\infty)$ and
$\rho\in
 (0,\gamma_1\wedge\gamma_2\wedge\frac{2-\beta}{2} )$,
there exists a positive constant $C$ such that
%
\begin{equation}
\label{s2.12}
\sup_{n\ge1} \bigl\llVert \bigl[X_n(t,x)-X_n(
\bar{t},\bar{x}) \bigr] 1_{L_n(\bar{t})} \bigr\rrVert _p \le C \bigl( |
\bar{t}- t |+ |\bar{x}-x | \bigr)^{\rho}.
\end{equation}
\end{theorem}

\begin{theorem}
\label{ts2.3}
Assume Hypothesis~\textup{\ref{H}}. Fix a compact set $K\subset\mathbb{R}^3$.
Then, for any
$p\in[1,\infty)$,
%
\begin{equation}
\label{s2.13}
\lim_{n\to\infty} \mathop{\sup_{t\in[0,T]}}_{x\in K(t)}\bigl\|
\bigl(X_n(t,x)-X(t,x) \bigr) 1_{L_n(t)}\bigr\|_p =0.
\end{equation}
\end{theorem}

The proof of Theorem~\ref{ts2.2} consists of two parts. First, we
shall consider $t=\bar t$ and obtain (\ref{s2.12}),
uniformly in $t\in[0,T]$. This is the difficult and novel part, and
where the additional assumption on $B$ being affine
is needed. Then, using this result, we consider $x=\bar x$ and
following the proof of \cite{Delgado--Sanz-Sole012}, Proposition~2.9,
we can establish (\ref{s2.12}), uniformly in $x$ over compact sets.
The details of the proof of the estimates of $L^p$
increments in time are omitted, since they can be reconstructed from
\cite{Delgado--Sanz-Sole012}, Proposition~2.9, with
minor changes.

\begin{remark}
\label{r2.1}
Assume that Hypothesis~\textup{\ref{H}} holds and, moreover,
\[
\sup_{n\ge1} \sup_{t\in[0,T]} \bigl\llVert
\bigl[X_n(t,x)-X_n(t,\bar x) \bigr]1_{L_n(t)} \bigr
\rrVert _p \le C |\bar{x}-x |^{\rho},
\]
with $\rho$ as in Theorem~\ref{ts2.2}. Then, with the same proof of
\cite{Delgado--Sanz-Sole012}, Proposition~2.9, one has
\[
\sup_{n\ge1} \bigl\llVert \bigl[X_n(t,x)-X_n(
\bar t,x) \bigr]1_{L_n(\bar t)} \bigr\rrVert _p \le C |\bar{t}-t
|^{\rho},
\]
for any $t_0\le t\le\bar{t}\le T$, $t_0>0$, uniformly over $x$ on a
compact set of $\mathbb{R}^3$.
\end{remark}

The proof of Theorem~\ref{ts2.3} is very similar to \cite{Delgado--Sanz-Sole012}, Theorem~2.4, and will also be omitted.
Notice that the initial condition $X^0(t,x)$ cancels in the difference
$X_n(t,x)-X(t,x)$, and also that in the proof
of \cite{Delgado--Sanz-Sole012}, Theorem~2.4, the stationarity property
is never used.


The rest of the section is devoted to establish $L^p$ estimates of
increments in space. They will be derived from
Proposition~\ref{pss2.1.2} below.

\begin{proposition}
\label{pss2.1.2}
We assume Hypothesis~\textup{\ref{H}} and that the function $B$ is affine. Fix a
compact set $K\subset\mathbb{R}^3$ and
$p\in (\frac{6\vee[2(4-\beta)]}{2-\beta}, \infty )$.
Then, for any $t\in[0,T]$ and
$\gamma\in (0,\gamma_1\wedge\gamma_2\wedge(\frac{2-\beta
}{2}-\frac{3}{p}) )$,
%
\begin{equation}
\label{s2.15} \sup_n\sup_{t\in[0,T]}
\mathbb{E} \bigl( \bigl[\bigl\Vert X_n(t)\bigr\Vert ^p_{W^{\gamma,p}(K(t))}
+ \bigl\Vert X_n^-(t)\bigr\Vert^p_{W^{\gamma,p}(K(t))}
\bigr]1_{L_n(t)} \bigr) < +\infty.
\end{equation}
\end{proposition}

Assume this has been proved.
By the Sobolev embedding theorem (see, for instance, \cite{Shimakura92}, Theorem E.12,\vspace*{1pt} page~257), for any bounded
or unbounded domain $\mathcal{O}\subset\mathbb{R}^d$, $W^{\rho
,p}(\mathcal
{O})\subset
\mathcal{C}^{\bar\rho}(\mathcal{O})$, for each $\bar\rho<\rho
-\frac{3}{p}$. Since Proposition~\ref{pss2.1.2}
holds for any $p$ large enough, (\ref{s2.15}) yields
%
\begin{equation}
\label{s2.14} \sup_{n\ge0} \sup_{t\in[0,T]}\bigl\|
\bigl(X_n(t,x)-X_n(t,\bar{x}) \bigr) 1_{L_n(t)}
\bigr\|_p \le C |x-\bar{x}|^{\rho},
\end{equation}
for any $\rho\in (0,\gamma_1\wedge\gamma_2\wedge\frac
{2-\beta}{2} )$.

The next Lemma~\ref{lss2.1.3} and Proposition~\ref{pss2.1.7} are
important ingredients in the proof of Proposition~\ref{pss2.1.2}.
Roughly speaking, Lemma~\ref{lss2.1.3} is an abstract
result about upper bounds of $L^p$ moments
of fractional norms of indefinite stochastic integral, taking into
account the size of the domain of integration
in time. In Proposition~\ref{pss2.1.7}, it is used to establish the
discrepancy in the fractional norm, and in terms
of $n$, between the Picard's iterations of $X_n(t,x)$ and
$X_n^{-}(t,x)$ (see (\ref{s2.4}), (\ref{s2.5}), resp.).

For a function $f \dvtx  \mathbb{R}^3 \mapsto\mathbb{R}$, we set
\begin{eqnarray*}
Df(u,x)&=& f(u+x)-f(u),
\nonumber
\\
D^2f(u,x)&= & f(u-x)-2f(u)+f(u+x).
\nonumber
\end{eqnarray*}

Given a bounded set $\mathcal{O}\in\mathbb{R}^3$ and $\varepsilon>0$, we
denote by $\mathcal{O}^\varepsilon$ the open set
%
\begin{equation}
\label{ext}
\mathcal{O}^\varepsilon= \bigl\{x\in\mathbb{R}^3
\dvt \exists z\in\mathcal{O} \mbox{ such that } |x-z|< \varepsilon
\bigr\}.
\end{equation}


\begin{lemma}
\label{lss2.1.3}
Fix $p\in (\frac{6}{2-\beta},\infty )$, $\gamma\in
 (0,\frac{2-\beta}{2}-\frac{3}{p} )$, $t\in(0,T]$
and a bounded domain $\mathcal{O}\subset R^3$.
Let $Z=\{Z(s,x), (s,x)\in[0,T]\times\mathbb{R}^3\}$ be a stochastic process
such that
%
\begin{equation}
\label{s2.42}
\int_0^t \mathrm{d}s \mathbb{E} \bigl(
\bigl\llVert Z(s) \bigr\rrVert _{W^{\gamma
,p}(\mathcal{O}^{s})}^p \bigr) <\infty.
\end{equation}
Then
%
\begin{eqnarray}
&& \int_{\mathcal{O}} \mathrm{d}x \int_{\mathcal{O}}
\mathrm{d}y \frac{\mathbb{E}
(\llVert  [G(\cdot,x-\ast)-G(\cdot,y-\ast) ] Z(\cdot,\ast
)\rrVert ^p_{\mathcal{H}_t} )}{\vert x-y\vert^{3+\gamma p}}
\nonumber
\\[-8pt]
\label{s2.4200}
\\[-8pt]
\nonumber
&&\quad \le C t^{\eta({p}/{2}-1)} \int_0^t
\mathrm{d}s \mathbb{E} \bigl( \bigl\llVert Z(s) \bigr\rrVert _{W^{\gamma,p}(\mathcal{O}^{s})}^p
\bigr),
\end{eqnarray}
with $\eta:=\inf (\frac{4-\beta}{2}, 3-2\gamma-\frac
{6}{p}-\beta )\in(1,2)$.
Consequently,
%
\begin{equation}
\label{s2.43}
\mathbb{E} \bigl( \bigl\llVert \bigl\llVert G(\cdot,\bullet-\ast)
Z(\cdot ,\ast) \bigr\rrVert _{\mathcal{H}_t} \bigr\rrVert _{\gamma,p,\mathcal
{O}}^p
\bigr)\le C t^{\eta({p}/{2}-1)} \int_0^t \mathrm{d}s
\mathbb{E} \bigl( \bigl\llVert Z(s) \bigr\rrVert _{W^{\gamma,p}(\mathcal{O}^{s})}^p
\bigr).
\end{equation}
\end{lemma}

\begin{pf}
Throughout the proof, $\beta\in(0,2)$ is fixed and we denote by
$f(x)$ the Riesz kernel $|x|^{-\beta}$. Remember that the symbols
``$\cdot$'', ``$\ast$'' denote the relevant variables for the $\mathcal{H}_t$
norm, and ``$\bullet$'' the argument for the fractional norm $\Vert
\cdot\Vert_{\gamma,p,\mathcal{O}}$.

Fix $x,y\in\mathbb{R}^3$. By applying the triangular inequality, we have
%
\begin{eqnarray}
&& \bigl\llvert \bigl\llVert G(\cdot,x-\ast) Z(\cdot,\ast) \bigr\rrVert
_{\mathcal{H}_t}- \bigl\llVert G(\cdot,y-\ast) Z(\cdot,\ast) \bigr\rrVert
_{\mathcal{H}_t} \bigr\rrvert
\nonumber
\\[-8pt]
\\[-8pt]
\nonumber
&&\quad\le \bigl\llVert \bigl[G(\cdot,x-\ast)- G(\cdot,y-
\ast) \bigr] Z(\cdot,\ast) \bigr\rrVert _{\mathcal{H}_t}.
\end{eqnarray}
Hence, (\ref{s2.43}) is a trivial consequence of (\ref{s2.4200}).

Set
\[
T_t:= \int_{\mathcal{O}} \mathrm{d}x \int_{\mathcal{O}}
\mathrm{d}y \frac{\mathbb
{E}
(\llVert  [G(\cdot,x-\ast)
-G(\cdot,y-\ast) ] Z(\cdot,\ast)\rrVert ^p_{\mathcal{H}_t
} )}{\vert x-y\vert^{3+\gamma p}}.
\]
By (\ref{mattila}), we write
%
\begin{eqnarray}
&& \bigl\llVert \bigl[G(\cdot,x-\ast)- G(\cdot,y-\ast) \bigr] Z(
\cdot,\ast) \bigr\rrVert _{\mathcal{H}_t}
\nonumber
\\
\label{s2.432}
&&\quad = C \biggl(\int_0^t
\!\int_{\mathbb{R}^3}\!\int_{\mathbb{R}^3} Z(s,u)f(u-v)Z(s,v)
 \\
\nonumber
&&\hspace*{56pt}\qquad{}\times\bigl[G(s,x-\mathrm{d}u)-G(s,y-\mathrm{d}u) \bigr]\bigl[G(s,x-\mathrm{d}v)-G(s,y-\mathrm{d}v) \bigr]
\biggr)^{{1}/{2}}.\qquad
\end{eqnarray}
Fix $p\in (\frac{6}{2-\beta},\infty )$, $\gamma\in
 (0,\frac{2-\beta}{2}-\frac{3}{p} )$
and $\mathcal{O}\subset\mathbb{R}^3$, and let $\bar\rho=\gamma
+\frac
{3}{p}$. By\vspace*{1pt} (\ref{s2.432}) and using the
method of the proof of \cite{Dalang--Sanz-Sole09}, Proposition~3.5,
increments of $G$ are transferred to
increments of the factors $f$ and $Z$. We obtain (see \cite{Dalang--Sanz-Sole09}, pages 19--20),
%
\begin{equation}
\label{s2.4321}
T_t\le C \sum_{i=1}^4
\mathbb{E} \biggl(\int_{\mathcal{O}}\mathrm{d}x \int_{\mathcal
{O}}\mathrm{d}y
\frac{\llvert  J_i^t(x,y)
\rrvert ^{{p}/{2}}}{|x-y|^{\bar\rho p}} \biggr),
\end{equation}
where
\begin{eqnarray*}
J_1^t(x,y)&=&\int_0^t
\mathrm{d}s \int_{\mathbb{R}^3}G(s,\mathrm{d}u)\int_{\mathbb{R}
^3}G(s,\mathrm{d}v)f(y-x+v-u)
\\
&&\hspace*{83pt}\quad{}\times
\bigl[Z(s,x-u)-Z(s,y-u) \bigr] \bigl[Z(s,x-v)-Z(s,y-v) \bigr],
\\
J_2^t(x,y)&=&\int_0^t
\mathrm{d}s \int_{\mathbb{R}^3}G(s,\mathrm{d}u)\int_{\mathbb{R}
^3}G(s,\mathrm{d}v)Df(v-u,x-y)Z(s,x-u)
\\
&&\hspace*{83pt}\quad{}\times
\bigl[Z(s,x-v)-Z(s,y-v) \bigr],
\\
J_3^t(x,y)&=&\int_0^t
\mathrm{d}s \int_{\mathbb{R}^3}G(s,\mathrm{d}u)\int_{\mathbb{R}
^3}G(s,\mathrm{d}v)Df(u-v,x-y)
\\
&&\hspace*{83pt}\quad{}\times
Z(s,x-v) \bigl[Z(s,x-u)-Z(s,y-u) \bigr],
\\
J_4^t(x,y)&=&\int_0^t
\mathrm{d}s \int_{\mathbb{R}^3}G(s,\mathrm{d}u)\int_{\mathbb{R}
^3}G(s,\mathrm{d}v)D^2f(v-u,x-y)
Z(s,x-u)Z(s,x-v).
\end{eqnarray*}

Let
%
\begin{equation}
\label{s2.433}
\mu_1(t,x,y):=\int_0^t
\mathrm{d}s \int_{\mathbb{R}^3}G(s,\mathrm{d}u)\int_{\mathbb
{R}^3}G(s,\mathrm{d}v)f(y-x+v-u).
\end{equation}
The following properties hold:
%
\begin{eqnarray} \label{s2.434}
\mu_1(t,x,y) &\le &  Ct^{3-\beta},
\\
\label{s2.435}
\sup_{s\in[0,T]}\int_{\mathbb{R}^3}G(s,\mathrm{d}u)\int
_{\mathbb{R}
^3}G(s,\mathrm{d}v)f(y-x+v-u) &\le &  C
\end{eqnarray}
(see, e.g., \cite{Dalang--Sanz-Sole13}, Lemma~5.1).
Thus, by H\"older's inequality and (\ref{s2.434}) we obtain
\begin{eqnarray*}
&&\mathbb{E} \biggl(\int_{\mathcal{O}}\mathrm{d}x \int_{\mathcal{O}}\mathrm{d}y
\frac
{\llvert  J_1^t(x,y)\rrvert ^{{p}/{2}}}{|x-y|^{\bar\rho
p}} \biggr)
\\
&&\quad\le C t^{(3-\beta)({p}/{2}-1)} \int_0^t \mathrm{d}s \int
_{\mathbb{R}
^3}G(s,\mathrm{d}u)\int_{\mathbb{R}^3}G(s,\mathrm{d}v)f(y-x+v-u)
\mathbb{E} \bigl(\bigl\Vert Z(s)\bigr\Vert ^p_{\gamma,p,\mathcal{O}^s} \bigr).
\end{eqnarray*}
By (\ref{s2.435}), this yields
%
\begin{equation}
\label{s2.437}
\mathbb{E} \biggl(\int_{\mathcal{O}}\mathrm{d}x \int
_{\mathcal{O}}\mathrm{d}y \frac
{\llvert  J_1^t(x,y)\rrvert ^{{p}/{2}}}{|x-y|^{\bar\rho
p}} \biggr) \le C t^{(3-\beta)({p}/{2}-1)}
\int_0^t \mathrm{d}s \bigl(\bigl\Vert Z(s)
\bigr\Vert^p_{\gamma,p,\mathcal{O}^s} \bigr).
\end{equation}

By symmetry, the contributions of the terms $J_2^t(x,y)$ and
$J_3^t(x,y)$ are the same. Hence, we will focus on $J_2^t(x,y)$.
Set
%
\begin{equation}
\label{s2.438} \mu_2(t,x,y)=\int_0^t
\mathrm{d}s \int_{\mathbb{R}^3}G(s,\mathrm{d}u)\int_{\mathbb{R}
^3}G(s,\mathrm{d}v)
\bigl\vert Df(v-u,x-y)\bigr\vert.
\end{equation}
The following properties hold:
%
\begin{eqnarray}\label{s2.439}
\mu_2(t,x,y) &\le &  C |x-y|^\alpha t^{3-(\alpha+\beta)},
\\
\label{s2.440}
\sup_{s\in[0,T]} \int_{\mathbb{R}^3}G(s,\mathrm{d}u)\int
_{\mathbb
{R}^3}G(s,\mathrm{d}v)\bigl\vert Df(v-u,x-y)\bigr\vert &\le &  C |x-y|^\alpha,
\end{eqnarray}
for any $\alpha\in(0,(2-\beta)\wedge1)$ (see \cite{Dalang--Sanz-Sole13}, Lemma~5.4, and a slight modification of
\cite{Dalang--Sanz-Sole09}, Lemma~6.1, resp.).
By applying H\"older's and Schwarz' inequalities and (\ref{s2.439}),
we can write
\begin{eqnarray*}
&&\!\!\! \mathbb{E} \biggl(\int_{\mathcal{O}}\mathrm{d}x \int_{\mathcal{O}}\mathrm{d}y
\frac
{\llvert  J_2^t(x,y)\rrvert ^{{p}/{2}}}{|x-y|^{\bar\rho
p}} \biggr)
\\
&&\!\!\!\quad=\mathbb{E} \biggl(\int_{\mathcal{O}}\mathrm{d}x \int_{\mathcal{O}}\mathrm{d}y
\biggl(\frac
{\llvert  J_2^t(x,y)\rrvert }{|x-y|^\alpha|x-y|^{2\bar\rho
-\alpha}} \biggr)^{{p}/{2}} \biggr)
\\
&&\!\!\!\quad\le\int_{\mathcal{O}}\mathrm{d}x \int_{\mathcal{O}}\mathrm{d}y \biggl(
\frac{\mu
_2(t,x,y)}{|x-y|^\alpha} \biggr)^{{p}/{2}-1} \int_0^t
\mathrm{d}s \int_{\mathbb{R}^3}G(s,\mathrm{d}u)
\\
&&\!\!\!\!\!\qquad{}\times\int_{\mathbb{R}^3}G(s,\mathrm{d}v)
\frac{\vert
Df(v-u,x-y)\vert
}{|x-y|^{\alpha}}\mathbb{E} \biggl(\bigl|Z(s,x-u)\bigr|^{{p}/{2}} \biggl(
\frac
{|Z(s,x-v)-Z(s,y-v)|}{|x-y|^{2\bar\rho-\alpha}} \biggr)^{{p}/{2}} \biggr)
\\
&&\!\!\!\quad\le C t^{[3-(\alpha+\beta)] [{p}/{2}-1 ]}
\\
&&\!\!\!\!\!\!{}\qquad\times\int_0^t \mathrm{d}s \biggl[\int
_{\mathcal{O}}\mathrm{d}x \int_{\mathcal{O}}\mathrm{d}y \int
_{\mathbb{R}^3}G(s,\mathrm{d}u)\int_{\mathbb{R}^3}G(s,\mathrm{d}v)
\frac{\vert
Df(v-u,x-y)\vert
}{|x-y|^{\alpha}}E \bigl(\bigl\vert Z(s,x-u)\bigr\vert^p \bigr)
\biggr]^{{1}/{2}}
\\
&&\!\!\!\!\!\!\qquad{}\times \biggl[\int_{\mathcal{O}}\mathrm{d}x \int_{\mathcal{O}}\mathrm{d}y
\int_{\mathbb{R}^3}G(s,\mathrm{d}u)
\\
&&\hspace*{12pt}\!\!\!\qquad{}\times\int_{\mathbb{R}^3}G(s,\mathrm{d}v)
\frac{\vert
Df(v-u,x-y)\vert
}{|x-y|^{\alpha}} E \biggl(\frac{\vert Z(s,x-v)-Z(s,y-v)\vert
}{|x-y|^{2\bar\rho-\alpha}} \biggr)^p
\biggr]^{{1}/{2}}.
\end{eqnarray*}
By applying (\ref{s2.440}), we conclude
%
\begin{eqnarray}
&& \mathbb{E} \biggl(\int_{\mathcal{O}}\mathrm{d}x \int
_{\mathcal{O}}\mathrm{d}y \frac
{\llvert  J_2^t(x,y)\rrvert ^{{p}/{2}}}{|x-y|^{\bar\rho p}} \biggr)
\nonumber
\\[-8pt]
\label{s2.441}
\\[-8pt]
\nonumber
&&\quad\le  C t^{[3-(\alpha+\beta)] [{p}/{2}-1 ]}\int_0^t \mathrm{d}s \bigl[\mathbb{E} \bigl(
\bigl\Vert Z(s)\bigr\Vert ^p_{L^p(\mathcal
{O}^s)} \bigr) \mathbb{E} \bigl(\bigl\Vert Z(s)
\bigr\Vert^p_{2\bar\rho-\alpha-{3}/{p},p,\mathcal{O}^s} \bigr) \bigr]^{{1}/{2}}.
\end{eqnarray}
Choose $\alpha=\gamma+\frac{3}{p}$. Remember that $\gamma<\frac
{2-\beta}{2}-\frac{3}{p}$. Hence
$\alpha\in (0,\frac{2-\beta}{2} )$. This implies $\alpha
\in(0,(2-\beta)\wedge1)$, as required.
Notice also that $2\bar\rho-\alpha-\frac{3}{p}=\gamma$, and
$3-(\alpha+\beta)>\frac{4-\beta}{2}>1$. Therefore,
%
\begin{eqnarray}
\label{s2.442}
&&\mathbb{E} \biggl(\int_{\mathcal{O}}\mathrm{d}x \int
_{\mathcal{O}}\mathrm{d}y \frac
{\llvert  J_2^t(x,y)\rrvert ^{{p}/{2}}}{|x-y|^{\bar\rho p}} \biggr) \le C t^{
(({4-\beta})/{2} ) ({p}/{2}-1 )}
\int_0^t \mathbb{E} \bigl(\bigl\Vert Z(s)
\bigr\Vert^p_{W^{\gamma,p}(\mathcal
{O}^s)} \bigr).
\end{eqnarray}

Set
\[
\mu_4(t,x,y)=\int_0^t \mathrm{d}s \int
_{\mathbb{R}^3}G(s,\mathrm{d}u)\int_{\mathbb{R}
^3}G(s,\mathrm{d}v)\bigl\vert
D^2f(v-u,x-y)\bigr\vert.
\]
The following properties hold:
%
\begin{eqnarray}\label{s2.444}
\mu_4(t,x,y) &\le &  C |x-y|^\alpha t^{3-(\alpha+\beta)},
\\
\label{s2.445}
\sup_{s\in[0,T]}\int_{\mathbb{R}^3}G(s,\mathrm{d}u)\int
_{\mathbb{R}
^3}G(s,\mathrm{d}v)D^2f(v-u,x-y) &\le &  C|x-y|^\alpha,
\end{eqnarray}
with $\alpha\in(0,2-\beta)$. The former is proved in
\cite{Dalang--Sanz-Sole13}, Lemma~5.5,  and the latter is a
slight modification of \cite{Dalang--Sanz-Sole09}, Lemma~6.2.

Choose $\alpha=2\bar\rho$. By applying H\"older's and Schwarz's
inequalities, we obtain
\begin{eqnarray*}
&&\!\!\! \mathbb{E} \biggl(\int_{\mathcal{O}}\mathrm{d}x \int_{\mathcal{O}}\mathrm{d}y
\frac
{\llvert  J_4^t(x,y)\rrvert ^{{p}/{2}}}{
|x-y|^{\bar\rho p}} \biggr)
\\
&&\!\!\!\quad = \int_{\mathcal{O}}\mathrm{d}x \int_{\mathcal{O}}\mathrm{d}y\mathbb{E}
\biggl(\frac{\llvert  J_4^t(x,y)\rrvert }{|x-y|^\alpha} \biggr)^{{p}/{2}}
\\
&&\!\!\!\quad \le C t^{[3-(\alpha+\beta)] [{p}/{2}-1 ]}
\\
&&\!\!\!\qquad{}\times\int_{\mathcal{O}}\mathrm{d}x \int_{\mathcal{O}}\mathrm{d}y \int
_0^t \mathrm{d}s \int_{\mathbb{R}^3}G(s,\mathrm{d}u)
\int_{\mathbb{R}^3}G(s,\mathrm{d}v)\frac{\vert
D^2f(v-u,x-y)\vert
}{|x-y|^\alpha}
\\
&&\!\!\!\hspace*{157pt}\qquad{}\times
E \bigl(\bigl|Z(s,x-u)\bigr|^{{p}/{2}}\bigl|Z(s,x-v)\bigr|^{{p}/{2}} \bigr)
\\
&&\!\!\!\quad\le C t^{[3-(\alpha+\beta)] [{p}/{2}-1 ]}
\\
&&\!\!\!\!\!\qquad{}\times\biggl(\int_{\mathcal{O}}\mathrm{d}x \int_{\mathcal{O}}\mathrm{d}y
\int_0^t \mathrm{d}s \int_{\mathbb{R}^3}G(s,\mathrm{d}u)
\int_{\mathbb{R}^3}G(s,\mathrm{d}v)\frac{\vert
D^2f(v-u,x-y)\vert
}{|x-y|^\alpha} \mathbb{E}
\bigl(\bigl|Z(s,x-u) \bigr|^p \bigr) \biggr)\!.
\end{eqnarray*}

The estimate (\ref{s2.445}) yields
%
\begin{eqnarray}
\mathbb{E} \biggl(\int_{\mathcal{O}}\mathrm{d}x \int_{\mathcal{O}}\mathrm{d}y
\frac
{\llvert  J_4^t(x,y)\rrvert ^{
{p}/{2}}}{|x-y|^{\bar\rho p}} \biggr)& \le& C t^{ (3-(2\bar\rho+\beta) ) ({p}/{2}-1 )} \int_0^t
\mathrm{d}s \mathbb{E} \bigl( \bigl\llVert Z(s) \bigr\rrVert ^p_{L^p(\mathcal
{O}^s)}
\bigr)
\nonumber
\\[-8pt]
\label{s2.446}
\\[-8pt]
\nonumber
&=&C t^{ (3-2\gamma-{6}/{p}-\beta ) ({p}/{2}-1 )} \int_0^t \mathrm{d}s \mathbb{E}
\bigl( \bigl\llVert Z(s) \bigr\rrVert ^p_{L^p(\mathcal
{O}^s)}
\bigr).
\end{eqnarray}

For\vspace*{1.5pt} $\gamma< \frac{2-\beta}{2}-\frac{3}{p}$, we have $\bar\eta:=
3-2\gamma-\frac{6}{p}-\beta>1$,
and for $\beta\in(0,2)$, $3-\beta>\frac{4-\beta}{2}$.
Hence, from the results proved so far, we see that (\ref{s2.4200})
follows from (\ref{s2.437}),
(\ref{s2.442}) and (\ref{s2.446}).
\end{pf}
%



For the proof of Proposition~\ref{pss2.1.2}, it is convenient to
consider localizations of the processes
$X_n$, $X_n^-$ in the space variable defined by $\{X_n(t,x)1_{K(t)}(x),
(t,x)\in[0,T]\times\mathbb{R}^3\}$,
$\{X_n^-(t,x)1_{K(t)}(x), (t,x)\in[0,T]\times\mathbb{R}^3\}$, respectively,
with $K(t)$ given in (\ref{loc}).

Let $x,y\in\mathbb{R}^3$ be such that $x\in K(t)$ and $|x-y|=t-s$. This
choice is motivated by the fact that the
Green function $G(t-s, x-\ast)$ has support on the sphere centered at
$x$ and with radius $t-s$. By the
triangular inequality, $d(y,K)\le d(y,x)+ d(x,K)\le T-s$. Thus, $y\in
K(s)$. Consequently,
$\{X_n(t,x)1_{K(t)}(x), (t,x)\in[0,T]\times\mathbb{R}^3\}$ satisfies the
following localized evolution
equation:
%
\begin{eqnarray}
X_n(t,x)1_{K(t)}(x) &=& X^0(t,x)1_{K(t)}(x)
\nonumber
\\
&&{}+1_{K(t)}(x)\int_0^t\! \int
_{\mathbb{R}^3} G(t-s,x-y) A \bigl(X_n(s,y)
\bigr)1_{K(s)}(y) M(\mathrm{d}s,\mathrm{d}y)
\nonumber
\\
\label{s2.16}
&&{}+1_{K(t)}(x) \bigl\langle G(t-\cdot,x-\ast)B \bigl(X_n(
\cdot,\ast) \bigr)1_{K(\cdot
)}(\ast),w^n \bigr
\rangle_{\mathcal{H}_t}
\\
&&{}+1_{K(t)}(x) \bigl\langle G(t-\cdot,x-\ast)D \bigl(X_n(
\cdot,\ast) \bigr)1_{K(\cdot
)}(\ast),h \bigr\rangle_{\mathcal{H}_t}
\nonumber
\\
&&{}+1_{K(t)}(x)\int_0^t \bigl[G(t-s,
\ast)\star b \bigl(X_n(s,\ast ) \bigr)1_{K(s)}(\ast)
\bigr](x) \,\mathrm{d}s.\nonumber
\end{eqnarray}
A similar equation also holds for $\{X_n^-(t,x)1_{K(t)}(x), (t,x)\in
[0,T]\times\mathbb{R}^3\}$, with the obvious changes.

Along with (\ref{s2.16}), we will also consider the Picard's
iterations defined by
%
\begin{eqnarray}
X_n^0(t,x)1_{K(t)}(x)&=&X^0(t,x)1_{K(t)}(x),
\nonumber
\\
X_n^m(t,x)1_{K(t)}(x) &=& X^0(t,x)1_{K(t)}(x)
\nonumber
\\
&&{}+1_{K(t)}(x)\int_0^t \!\int
_{\mathbb{R}^3} G(t-s,x-y) A \bigl(X_n^{m-1}(s,y)
\bigr)1_{K(s)}(y) M(\mathrm{d}s,\mathrm{d}y)
\nonumber
\\
\label{s2.17}
&&{}+1_{K(t)}(x) \bigl\langle G(t-\cdot,x-\ast)B \bigl(X_n^{m-1}(
\cdot,\ast ) \bigr)1_{K(\cdot)}(\ast),w^n \bigr
\rangle_{\mathcal{H}_t}
\\
&&{}+1_{K(t)}(x) \bigl\langle G(t-\cdot,x-\ast)D \bigl(X_n^{m-1}(
\cdot,\ast ) \bigr)1_{K(\cdot)}(\ast),h \bigr\rangle_{\mathcal{H}_t}
\nonumber
\\
\nonumber
&&{}+1_{K(t)}(x)\int_0^t \bigl[G(t-s,
\ast)\star b \bigl(X_n^{m-1}(s,\ast ) \bigr)1_{K(s)}(
\ast) \bigr](x) \,\mathrm{d}s,\qquad m\ge1.
\end{eqnarray}

For these Picard's iterations, and similarly as in (\ref{s2.5}), we define
%
\begin{eqnarray}
&& X_n^{-,0}(t,x)1_{K(t)}(x)=X^0(t,x)1_{K(t)}(x),
\nonumber
\\
&& X_n^{-,m}(t,x)1_{K(t)}(x)\nonumber\\
&&\quad=X^0(t,x)1_{K(t)}(x)
\nonumber
\\
\label{s2.170}
&&\qquad{}+1_{K(t)}(x)\int_0^{t_n}\!\! \int
_{\mathbb{R}^3} G(t-s,x-y) A \bigl(X_n^{m-1}(s,y)
\bigr)1_{K(s)}(y) M(\mathrm{d}s,\mathrm{d}y)
\\
\nonumber
&&\qquad{}+1_{K(t)}(x) \bigl\langle G(t-\cdot,x-\ast)B \bigl(X_n^{m-1}(
\cdot,\ast ) \bigr)1_{K(\cdot)}(\ast) 1_{[0,t_n]}(
\cdot),w^n \bigr\rangle_{\mathcal
{H}_t}
\\
&&\qquad{}+1_{K(t)}(x) \bigl\langle G(t-\cdot,x-\ast)D \bigl(X_n^{m-1}(
\cdot,\ast ) \bigr)1_{K(\cdot)}(\ast) 1_{[0,t_n]}(\cdot),h \bigr
\rangle_{\mathcal{H}_t}
\nonumber
\\
&&\qquad{}+1_{K(t)}(x)\int_0^{t_n} \bigl[G(t-s,
\ast)\star b \bigl(X_n^{m-1}(s,\ast ) \bigr)1_{K(s)}(
\ast) \bigr](x) \,\mathrm{d}s,\qquad m\ge1.\nonumber
\end{eqnarray}
%

In the next proposition, we consider the stochastic processes $\{X_n^m
(t,x), (t,x)\in[0,T]\times\mathbb{R}^3\}$,
$\{X_n^{-,m} (t,x), (t,x)\in[0,T]\times\mathbb{R}^3\}$ given in
(\ref{s2.17}), (\ref{s2.170}), respectively.

\begin{proposition}
\label{pss2.1.7}
Let $p>\frac{6\vee[2(4-\beta)]}{2-\beta}$ and $\gamma$ be as in
Lemma~\ref{lss2.1.3}. We also assume that the function $B$ is affine.
Fix $m\ge1$ and assume that
%
\begin{equation}
\label{s2.480}
\sup_{t\in[0,T]}\mathbb{E} \bigl( \bigl[\bigl\Vert
X_n^{m-1}(t)\bigr\Vert ^p_{W^{\gamma,p}(K(t))} + \bigl\Vert
X_n^{-,m-1}(t) \bigr\Vert^p_{W^{\gamma,p}(K(t))}
\bigr]1_{L_n(t)} \bigr) \le C,
\end{equation}
for some constant $C$ independent of $n$, $m$.

Then there exists $\bar\eta\in(1,\infty)$ independent of $n$, $m$
but depending on $p$, and $C>0$, such that
%
\begin{equation}
\label{s2.48}
\sup_{t\in[0,T]} \mathbb{E} \bigl( \bigl\llVert
X_n^{m}(t)-X_n^{-,m}(t) \bigr\rrVert
^p_{\gamma,p,K(t)}1_{L_n(t)} \bigr) \le C 2^{-n\bar\eta{p}/{2}}.
\end{equation}
\end{proposition}

\begin{pf}
Fix\vspace*{1pt} $p>\frac{6}{2-\beta}$, $\gamma\in (0,\frac{2-\beta
}{2}-\frac{3}{p} )$, $m\ge1$, $t\in(0,T]$. Remember that
if $x\in K(t)$ and $|x-y|=t-s$, then $y\in K(s)$. Hence, from (\ref
{s2.17}), (\ref{s2.170}), we have
%
\begin{equation}
\label{s2.481}
\mathbb{E} \bigl( \bigl\llVert X_n^{m}(t)-X_n^{-,m}(t)
\bigr\rrVert ^p_{\gamma
,p,K(t)}1_{L_n(s)} \bigr)\le C \sum
_{i=1}^5 V_{n,m}^i(t),
\end{equation}
where
\begin{eqnarray*}
V_{n,m}^1(t)&=&\mathbb{E} \biggl( \biggl\llVert \int
_{t_n}^t\!\int_{\mathbb{R}^3} G(t-s,
\bullet-y)A \bigl(X_n^{m-1}(s,y) \bigr) M(\mathrm{d}s,\mathrm{d}y) \biggr
\rrVert ^p_{\gamma,p,K(t)}1_{L_n(t)} \biggr),
\\
V_{n,m}^2(t)&=&\mathbb{E} \bigl( \bigl\llVert \bigl\langle
G(t-\cdot ,\bullet -\ast)B \bigl(X_n^{-,m-1}(\cdot,\ast)
\bigr) 1_{[t_n,t]}(\cdot),w^n \bigr\rangle_{\mathcal{H}_t} \bigr
\rrVert ^p_{\gamma
,p,K(t)}1_{L_n(t)} \bigr),
\\
V_{n,m}^3(t)&=&\mathbb{E} \bigl(\bigl\Vert \bigl\langle G(t-
\cdot, \bullet -\ast)
\bigl[B \bigl(X_n^{-,m-1}(\cdot,\ast)
\bigr) - B \bigl(X_n^{m-1}(\cdot,\ast) \bigr) \bigr]
\\
&&\hspace*{167pt}{}\times
1_{[t_n,t]}(\cdot),w^n \bigr\rangle_{\mathcal{H}_t}
\bigr\Vert^p_{\gamma,p,K(t)}1_{L_n(t)} \bigr),
\\
V_{n,m}^4(t)&=&\mathbb{E} \bigl( \bigl\llVert \bigl\langle
G(t-\cdot ,\bullet -\ast) D \bigl(X_n^{m-1}(\cdot,\ast)
\bigr)1_{[t_n,t]}(\cdot),h \bigr\rangle _{\mathcal{H}_t
} \bigr\rrVert
^p_{\gamma,p,K(t)}1_{L_n(t)} \bigr),
\\
V_{n,m}^5(t)&=&\mathbb{E} \biggl( \biggl\llVert \int
_{t_n}^t \bigl[G(t-s,\ast )\star b
\bigl(X_n^{m-1} (s,\ast) \bigr) \bigr](\bullet) \,\mathrm{d}s \biggr
\rrVert ^p_{\gamma
,p,K(t)}1_{L_n(t)} \biggr).
\end{eqnarray*}

Set\vspace*{1.5pt} $\bar\rho=\gamma+\frac{3}{p}$. By writing explicitly the norm
$\Vert\cdot\Vert_{\gamma,p,K(t)}$, and then applying Fubini's
theorem and Burkholder's inequality, we obtain
\begin{eqnarray*}
V_{n,m}^1(t) &\le &  C \int_{K(t)} \mathrm{d}x
\\
&&{}\times\int
_{K(t)} \mathrm{d}z\frac{\mathbb{E} (\llVert [G(t-\cdot,x-\ast
)-G(t-\cdot,z-\ast)] A(X_n^{m-1}(\cdot,\ast))1_{[t_n,t]}(\cdot
)1_{L_n(t)}\rrVert _{\mathcal{H}_t}^p )}{|x-z|^{\bar\rho p}}.
\end{eqnarray*}
Set $Z(s,y):=A(X_n^{m-1}(t-s,y))1_{L_n(t-s)}$. With the change of
variable $s\mapsto t-s$, the preceding inequality implies
\[
V_{n,m}^1(t)\le C \int_{K(t)} \mathrm{d}x \int
_{K(t)} \mathrm{d}z \frac{\mathbb
{E} (\int_0^{(2T2^{-n})\wedge t}
\,\mathrm{d}s (\llVert [G(s,x-\ast)-G(s,z-\ast)] Z(\cdot,\ast)\rrVert
^2_\mathcal{H} )^{{p}/{2}}}{|x-z|^{\bar\rho p}}. 
\]
The right-hand side of the preceding expression coincides up to a
constant with the left-hand side of (\ref{s2.4200}) with
$t:=(2T2^{-n})\wedge t$ and $\mathcal{O}:=K(t)$.

We\vspace*{1pt} are assuming
$\sup_{t\in[0,T]}\mathbb{E} (\llVert  X^{m-1}_n(t)\rrVert
^p_{W^{\gamma,p}(K(t))} 1_{L_n(t)} )<\infty$.
Hence, using the linear growth of $A$ and Lemma~\ref{lss2.1.1} [see
(\ref{s2.39})], we have
\begin{eqnarray*}
&& \sup_{s\in[0,t]} \mathbb{E} \bigl( \bigl\llVert A
\bigl(X_n^{m-1}(t-s) \bigr)1_{L_n(t-s)} \bigr\rrVert
^p_{W^{\gamma
,p}((K(t))^s)} \bigr)
\\
&&\quad\le\sup_{s\in[0,t]} C \bigl(\mathbb{E} \bigl( \bigl\llVert A
\bigl(X_n^{m-1}(t-s) \bigr) 1_{L_n(t-s)} \bigr\rrVert
^p_{L^p((K(t))^s)} \bigr)
\\
&&\hspace*{34pt}\qquad
{}+ \mathbb{E} \bigl( \bigl\llVert A \bigl(X_n^{m-1}(t-s)
\bigr)1_{L_n(t-s)} \bigr\rrVert ^p_{p,\gamma,(K(t))^s} \bigr) \bigr)
\\
&&\quad\le C_1+C_2 \sup_{s\in[0,t]} \mathbb{E}
\bigl( \bigl\llVert A \bigl(X_n^{m-1}(s)
\bigr)1_{L_n(s)} \bigr\rrVert ^p_{p,\gamma,(K(t))^{t-s}} \bigr)
\\
&&\quad\le C_1+C_2 \sup_{s\in[0,t]} \mathbb{E}
\bigl( \bigl\llVert X_n^{m-1}(s) 1_{L_n(s)} \bigr
\rrVert ^p_{\gamma,p,K(s)} \bigr)
\\
&&\quad \le C.
\end{eqnarray*}
%
By Lemma~\ref{lss2.1.3}, we conclude
%
\begin{equation}
\label{s2.482}
V_{n,m}^1(t)\le C 2^{-({np}/{2})\eta},
\end{equation}
with $\eta=\inf (\frac{4-\beta}{2}, 3-2\gamma-\frac
{6}{p}-\beta )$.


The term $V_{n,m}^2(t)$ is also a stochastic integral with respect to
$M$. Indeed,
for a given function $f \dvtx  [0,T]\times\mathbb{R}^3\to\mathbb{R}$ and
$t\in[0,T]$, let
$\tau_n$ be the operator defined by
%
\begin{equation}
\label{s3.34} \tau_n(f) (s,x)=f \bigl( \bigl(s+2^{-n}
\bigr)\wedge t, x \bigr).
\end{equation}
Let $\mathcal{E}_n$ be the closed subspace of $\mathcal{H}_T$
generated by the orthonormal system
\[
2^nT^{-1}1_{\Delta_i}(\cdot)\otimes e_j(
\ast),\qquad i=0,\ldots ,2^n-1, j=1,\ldots,n,
\]
and denote by $\pi_n$ the orthogonal projection operator on $\mathcal
{E}_n$. Notice that $\pi_n\circ\tau_n$ is
a bounded operator on $\mathcal{H}_T$, uniformly in $n$.

The random vector $X_n^{-,m-1}(s,\ast)$ is $\mathcal
{F}_{s_n}$-measurable. Then, using the definition of $w^n$,
it is easy to check that
\begin{eqnarray*}
&& \bigl\langle G(t-\cdot,\bullet-\ast)B \bigl(X_n^{-,m-1}(
\cdot,\ast ) \bigr)1_{[t_n,t]}(\cdot),w^n \bigr
\rangle_{\mathcal{H}_t}
\\[-1pt]
& &\quad= \int_{t_n}^t \!\int_{\mathbb{R}^3}
(\pi_n\circ\tau_n) \bigl(G(t-\cdot,\bullet-\ast) B
\bigl(X_n^{-,m-1}(\cdot,\ast) \bigr) \bigr) (s,y) M(\mathrm{d}s,\mathrm{d}y).
\end{eqnarray*}
Therefore, $V_{n,m}^2(t)$ can be studied in a similar way than
$V_{n,m}^1(t)$, with
$Z(s,y):=B(X_n^{-,m-1}(t-s,y))1_{L_n(t-s)}$. We obtain
%
\begin{equation}
\label{v2} V_{n,m}^2(t)\le C 2^{-({np}/{2})\eta},
\end{equation}
with $\eta=\inf (\frac{4-\beta}{2}, 3-2\gamma-\frac
{6}{p}-\beta )$.

Using Schwarz's inequality and (\ref{s2.10}), we have
%
\begin{eqnarray*}
V_{n,m}^3(t) &\le &   C n^{{3p}/{2}} 2^{{np}/{2}}
\mathbb{E} \bigl( \bigl\Vert \bigl\Vert G(t-\cdot,\bullet-\ast) \bigl[B
\bigl(X_n^{-,m-1} (\cdot,\ast) \bigr)- B \bigl(X_n^{m-1}(
\cdot,\ast) \bigr) \bigr]
\\[-1pt]
&&\hspace*{135pt}\qquad\qquad{}\times
1_{[t_n,t]}(\cdot) \bigr\Vert_{\mathcal{H}_t} 1_{L_n(t)}
\bigr\Vert_{\gamma,p,K(t)}^p \bigr).
\end{eqnarray*}
%
Let
%
\begin{equation}
\label{s2.4280}
Z(s,y):= \bigl[B \bigl(X_n^{-,m-1}(t-s,y)
\bigr)- B \bigl(X_n^{m-1}(t-s,y) \bigr) \bigr]1_{L_n(t-s)}.
\end{equation}
With the change of variable $s\mapsto t-s$, we see that
\[
V_{n,m}^3(t)\le C n^{{3p}/{2}} 2^{{np}/{2}}
\mathbb{E} \bigl( \bigl\llVert \bigl\llVert G(\cdot,\bullet-\ast) Z(\cdot,\ast)
\bigr\rrVert _{\mathcal{H}_{(2T2^{-n})\wedge t}} \bigr\rrVert ^p_{\gamma,p,K(t)} \bigr).
\]
We are assuming (\ref{s2.480}).
Hence, as in the analysis of $V_{n,m}^1(t)$, using Lemma~\ref
{lss2.1.1} we can prove that the assumptions of
Lemma~\ref{lss2.1.3} are satisfied for $\mathcal{O}:=K(t)$. Consequently,
%
\begin{equation}
\label{s2.483}
V_{n,m}^3(t)\le C n^{{3p}/{2}}
2^{{np}/{2}} 2^{-n\eta({p}/{2}-1)} \int_0^{(2T2^{-n})\wedge t} E
\bigl( \bigl\llVert Z(s) \bigr\rrVert ^p_{W^{\gamma,p}(K(t)^s} \bigr),
\end{equation}
with $Z$ given in (\ref{s2.4280}) and $\eta=\inf (\frac
{4-\beta}{2}, 3-2\gamma-\frac{6}{p}-\beta )$.

From (\ref{ext}), it follows that for $0\le s\le t\le T$,
$K(t)^s=K(t-s)$. Therefore,
\begin{eqnarray*}
&& E \bigl( \bigl\llVert \bigl[B \bigl(X_n^{-,m-1}(t-s,y)
\bigr)- B \bigl(X_n^{m-1}(t-s,y) \bigr) \bigr]1_{L_n(t-s)}
\bigr\rrVert _{L^p((K(t)^s))}^p \bigr)
\\[-1pt]
&&\quad \le C \int_{K(t-s)} \mathrm{d}y E \bigl( \bigl\llvert
X_n^{-,m-1}(t-s,y)-X_n^{m-1}(t-s,y) \bigr
\rrvert ^p1_{L_n(t-s)} \bigr)
\\[-1pt]
&&\quad \le C \sup_{(s,y)\in[0,T]\times\mathbb{R}^3} E \bigl( \bigl\llvert X_n^{-,m-1}(s,y)-X_n^{m-1}(s,y)
\bigr\rrvert ^p 1_{L_n(s)} \bigr)
\\[-1pt]
&&\quad \le C n^{{3p}/{2}} 2^{-np({3-\beta})/{2}},
\end{eqnarray*}
where\vadjust{\goodbreak} we have applied Lemma~\ref{lss2.1.4}.

Since $B$ is affine, we have
\begin{eqnarray*}
&& \bigl\llVert \bigl[B \bigl(X_n^{-,m-1}(s) \bigr)- B
\bigl(X_n^{m-1}(s) \bigr) \bigr]1_{L_n(s)} \bigr
\rrVert ^p_{\gamma,p,K(s)}
\\
&&\quad \le C \bigl\llVert \bigl[X_n^{-,m-1}(s)-
X_n^{m-1}(s) \bigr]1_{L_n(s)} \bigr\rrVert
^p_{\gamma,p,K(s)},\qquad s\in[0,T].
\end{eqnarray*}
Applying these estimates to (\ref{s2.483}) yields
%
\begin{eqnarray}
V_{n,m}^3(t) & \le& C n^{{3p}/{2}}
2^{{np}/{2}} 2^{-n\eta
({p}/{2}-1)}\nonumber \\
&&{}\times\biggl[n^{{3p}/{2}} 2^{-np({3-\beta})/{2}}
 +\int_0^{(2T2^{-n})\wedge t} \mathrm{d}s \mathbb{E} \bigl( \bigl\llVert
X_n^{-,m-1}(s)-X_n^{m-1}(s) \bigr\rrVert
_{\gamma,p,K(s)}^p \bigr) \biggr]
\nonumber
\\[-8pt]
\label{s2.484}
\\[-8pt]
\nonumber
&\le &  C_1 n^{3p}2^{-({np}/{2}) (\eta+2-\beta-{2\eta}/{p} )}
\\
\nonumber
&&{}+ C_2 n^{{3p}/{2}}2^{-({np}/{2}) (\eta-1-{2\eta}/{p } )}\int
_0^{(2T2^{-n})\wedge t} \mathrm{d}s \mathbb{E}
\bigl( \bigl\llVert X_n^{-,m-1}(s)-X_n^{m-1}(s)
\bigr\rrVert _{\gamma,p,K(s)}^p \bigr).
\end{eqnarray}

Schwarz' inequality implies
\[
V_{n,m}^4(t)\le C \mathbb{E} \bigl( \bigl\llVert \bigl
\llVert G(\cdot ,\bullet -\ast) Z(\cdot,\ast) \bigr\rrVert _{\mathcal{H}_{(2T2^{-n})\wedge t}} \bigr
\rrVert ^p_{\gamma,p,K(t)} \bigr),
\]
with $Z(s,y):=D(X_n^{m-1}(t-s,y))1_{L_n(t-s)}$. Therefore, similarly as
in the study of the term $V_{n,m}^1(t)$ we obtain
%
\begin{equation}
\label{s2.485}
V_{n,m}^4(t)\le C 2^{-({np}/{2})\eta},
\end{equation}
with $\eta=\inf (\frac{4-\beta}{2}, 3-2\gamma-\frac
{6}{p}-\beta )$.

To study $V_{n,m}^5(t)$, we first apply Minkowski's inequality and then
the linear growth of $b$ and Lemma~\ref{lss2.1.1}.
We obtain (the details are left to the reader),
%
\begin{equation}
\label{s2.486}
V_{n,m}^5(t)\le C 2^{-np}.
\end{equation}

With (\ref{s2.481}), (\ref{s2.482}), (\ref{v2}), (\ref{s2.484}),
(\ref{s2.485}), (\ref{s2.486}), we have
\begin{eqnarray*}
&& \mathbb{E}  \bigl( \bigl\llVert X_n^m(t)-X^{-,m}(t)
\bigr\rrVert _{\gamma
,p,K(t)}^p \bigr)\\
&&\quad \le c_1
2^{-({np}/{2})\eta} +c_2 n^{3p}2^{-({np}/{2}) (\eta
+2-\beta-{2\eta}/{p} )}
\\
&&\qquad {}+c_3n^{{3p}/{2}}2^{-({np}/{2}) (\eta-1-{2\eta}/{p } )}\int_0^{(2T2^{-n})\wedge t}
\mathrm{d}s \mathbb{E} \bigl( \bigl\llVert X_n^{-,m-1}(s)-X_n^{m-1}(s)
\bigr\rrVert _{\gamma,p,K(s)}^p \bigr),
\end{eqnarray*}
with $\eta=\inf (\frac{4-\beta}{2}, 3-2\gamma-\frac
{6}{p}-\beta )\in(1,2)$ (see Lemma~\ref{lss2.1.3}).

In Lemma~\ref{aux1}, we prove that
for any $p>\frac{2(4-\beta)}{2-\beta}$, $\eta_1:=\eta-1-\frac
{2\eta}{p}>0$. This implies $\eta+2-\beta-\frac{2\eta}{p}>1$.

Set $f_n^m(t)=\mathbb{E} (\llVert  X_n^m(t)-X^{-,m}(t)\rrVert
_{\gamma,p,K(t)}^p )$, and $\eta_2:=\inf (\eta,\eta
+2-\beta-\frac{2\eta}{p} )$. We have proved that
\[
f_n^m(t)\le\varphi_n+\psi_n\int
_0^{2T2^{-n}\wedge t} \mathrm{d}s f_n^m(s),
\]
with $\varphi_n:=(c_1\vee c_2) 2^{-({np}/{2})\eta_2}$, $\psi
_n:=c_3n^{{3p}/{2}}2^{-np\eta_1}$, $\eta_1>0$, $\eta_2>1$. Then,
Gronwall's lemma yields
\[
f_n^m(t)\le\varphi_n \bigl(1+\exp(T
\psi_n) \bigr).
\]
Clearly, $\sup_n \psi_n<\infty$. Thus,
\[
f_n^m(t)\le C \varphi_n,\qquad n\ge1,
\]
and this yields (\ref{s2.48}).
\end{pf}
%


\begin{pf*}{Proof of Proposition~\protect\ref{pss2.1.2}}
Fix $p$ and $\gamma$ as in the assertion.
Using induction on $m\ge0$, we will first establish a result analogue
to (\ref{s2.15}) for the Picard's
iterations defined in (\ref{s2.17}) and (\ref{s2.170}). More
precisely, we will prove
%
\begin{equation}
\label{s2.19}
\sup_{t\in[0,T]}\mathbb{E} \bigl( \bigl[\bigl\Vert
X_n^m(t)\bigr\Vert ^p_{W^{\gamma
,p}(K(t))} + \bigl\Vert
X_n^{-,m}(t)\bigr\Vert^p_{W^{\gamma,p}(K(t))}
\bigr]1_{L_n(t)} \bigr) \le C,
\end{equation}
for some constant $C$ independent of $n,m$. By Fatou's lemma, and the
convergence in the $L^p$ norm of $X_n^m(t,x) 1_{L_n(t)}$,
and $X_n^{-,m}(t,x) 1_{L_n(t)}$ to $X_n(t,x) 1_{L_n(t)}$,
and $X_n^{-}(t,x) 1_{L_n(t)}$, respectively, for any $(t,x)\in
[0,T]\times\mathbb{R}^3$ [see (\ref{s2.410})],
this will imply (\ref{s2.15}).

For $m=0$, (\ref{s2.19}) is just the property (H2.1), which is a
consequence of hypothesis~(H2).

Let $m\ge1$ and assume that (\ref{s2.19}) holds for any Picard
iterative of order less or equal than $m-1$.
We recall that if $x\in K(t)$ and $|x-y|=t-s$, then $y\in K(s)$. Thus,
from (\ref{s2.17}) we see that
\[
\mathbb{E} \bigl(\bigl\Vert X_n^m(t)\bigr\Vert^p_{W^{\gamma
,p}(K(t))}1_{L_n(t)}
\bigr) \le C\sum_{i=1}^6
R_{n,m}^i(t),
\]
with
\begin{eqnarray*}
R_{n,m}^1(t)&=& \bigl\Vert X^0(t)
\bigr\Vert^p_{W^{\gamma,p}(K(t))},
\\
R_{n,m}^2(t)&=& \mathbb{E} \biggl( \biggl\llVert \int
_0^t \!\int_{\mathbb{R}^3} G(t-s,
\bullet-y) A \bigl(X_n^{m-1}(s,y) \bigr) M(\mathrm{d}s,\mathrm{d}y) \biggr
\rrVert ^p_{W^{\gamma
,p}(K(t))}1_{L_n(t)} \biggr),
\\
R_{n,m}^3(t)&=& \mathbb{E} \bigl( \bigl\llVert \bigl\langle
G(t-\cdot,\bullet -\ast )B \bigl(X_n^{-,m-1}(\cdot,\ast)
\bigr),w^n \bigr\rangle_{\mathcal{H}_t} \bigr\rrVert
^p_{W^{\gamma,p}(K(t))}1_{L_n(t)} \bigr),
\\
R_{n,m}^4(t)&=& \mathbb{E} \bigl(\bigl\Vert \bigl\langle G(t-
\cdot,\bullet-\ast ) \bigl[B \bigl(X_n^{m-1}(\cdot,\ast)
\bigr)
 -B \bigl(X_n^{-,m-1}(\cdot,\ast) \bigr)
\bigr],w^n \bigr\rangle _{\mathcal{H}_t}\bigr\Vert^p_{W^{\gamma,p}(K(t))}1_{L_n(t)}
\bigr),
\\
R_{n,m}^5(t)&=& \mathbb{E} \bigl( \bigl\llVert \bigl
\langle G(t-\cdot,\bullet -\ast )D \bigl(X_n^{m-1}(\cdot,
\ast) \bigr),h \bigr\rangle_{\mathcal{H}_t} \bigr\rrVert ^p_{W^{\gamma
,p}(K(t))}1_{L_n(t)}
\bigr),
\\
R_{n,m}^6(t)&=& \mathbb{E} \biggl( \biggl\llVert \int
_0^t \bigl[G(t-s,\ast )\star b
\bigl(X_n^{m-1}(s,\ast) \bigr) \bigr](\bullet) \,\mathrm{d}s \biggr
\rrVert ^p_{W^{\gamma
,p}(K(t))}1_{L_n(t)} \biggr),
\end{eqnarray*}
where the symbols ``$\cdot$'', ``$\ast$'' denote the time and space
variables, respectively, that are relevant
for the $\mathcal{H}_t$ norm, while the symbol ``$\bullet$'' denotes
the argument corresponding to functions in
the space $W^{\gamma,p}(K(t))$.

As has been pointed out before, the assumption (H2) implies
%
\begin{equation}
\label{s2.18} R_{n,m}^1(t)\le C.
\end{equation}

By the induction hypotheses and (\ref{s2.38}), (\ref{s2.39}) applied
to the function $g:=A$ and $Z(t,x):=X_n^{m-1}(t,x) 1_{L_n(t)}$, we see
that the assumptions of \cite{Dalang--Sanz-Sole09}, Theorem~3.1, hold.
Therefore, we have
\[
R_{n,m}^2(t)\le C \int_0^t
\mathrm{d}s \mathbb{E} \bigl( \bigl\llVert A \bigl(X_n^{m-1}(s)
\bigr) \bigr\rrVert ^p_{W^{\gamma
,p}(K(t)^{t-s})}1_{L_n(s)} \bigr).
\]
From this inequality, the definition (\ref{sobolevfrac}), the
Lipschitz continuity of the function $A$ and Lemma~\ref{lss2.1.1}, it
follows that
%
\begin{equation}
\label{s2.20} R_{n,m}^2(t)\le C_1+C_2
\int_0^t \mathrm{d}s \mathbb{E} \bigl( \bigl\llVert
X_n^{m-1}(s) \bigr\rrVert ^p_{W^{\gamma,p}(K(s))}1_{L_n(s)}
\bigr)\le C.
\end{equation}

As has been mentioned in the analysis of the term $V_{n,m}^2(t)$ in
Proposition~\ref{pss2.1.7},
the following identity holds:
\begin{eqnarray*}
&& \bigl\langle G(t-\cdot,\bullet-\ast)B \bigl(X_n^{-,m-1}(
\cdot,\ast ) \bigr),w^n \bigr\rangle_{\mathcal{H}_t}
\\
&&\quad= \int_0^t \!\int_{\mathbb{R}^3} (
\pi_n\circ\tau_n) (G(t-\cdot ,\bullet-\ast) B
\bigl(X_n^{-,m-1}(\cdot,\ast) \bigr) (s,y) M(\mathrm{d}s,\mathrm{d}y).
\end{eqnarray*}
Consequently, $R_{n,m}^3(t) \le C
[R_{n,m}^{3,1}(t)+R_{n,m}^{3,2}(t) ]$,
with
\begin{eqnarray*}
&& R_{n,m}^{3,1}(t)
\\
&&\quad =\mathbb{E} \biggl( \biggl\llVert \int_0^t
\!\int_{\mathbb{R}^3}(\pi _n\circ\tau _n)
\bigl(G(t-\cdot,\bullet-\ast) B \bigl(X_n^{-,m-1}(\cdot,\ast)
\bigr) \bigr) (s,y) M(\mathrm{d}s,\mathrm{d}y) \biggr\rrVert _{L^p(K(t))}^p1_{L_n(t)}
\biggr),
\\
&&R_{n,m}^{3,2}(t)
\\
&&\quad =\mathbb{E} \biggl( \biggl\llVert \int_0^t\!
\int_{\mathbb{R}^3}(\pi _n\circ\tau _n)
\bigl(G(t-\cdot,\bullet-\ast) B \bigl(X_n^{-,m-1}(\cdot,\ast)
\bigr) \bigr) (s,y) M(\mathrm{d}s,\mathrm{d}y) \biggr\rrVert _{\gamma,p,K(t)}^p1_{L_n(t)}
\biggr).
\end{eqnarray*}

By developing the $L^p(K(t))$ norm, and using Fubini's theorem,
Burkholder's inequality and the boundedness
of the operator $\pi_n\circ\tau_n$, we have
%
\begin{equation}
\label{r31}
R_{n,m}^{3,1}(t) \le C \mathbb{E} \biggl(
\biggl(\int_{K(t)} \mathrm{d}x \bigl\llVert \bigl(G(t-\cdot,x-\ast)B
\bigl(X_n^{-,m-1} \bigr)  (\cdot,\ast) \bigr) \bigr
\rrVert _{\mathcal{H}_t}^p \biggr)1_{L_n(t)} \biggr).
\end{equation}
Now we apply the usual estimates on the $\mathcal{H}_t$ norm along
with the
property (\ref{s2.38}) and the
induction assumption to conclude that
$R_{n,m}^{3,1}(t) \le C$.

Using the definition of the fractional norm (see (\ref{frac})),
Fubini's theorem, along with
Burkholder's inequality and the boundedness of the operator $\pi
_n\circ\tau_n$ yields,
%
\begin{eqnarray}
R_{n,m}^{3,2}(t)& =& \int_{K(t)} \mathrm{d}x
\int_{K(t)}\mathrm{d}z \frac
{1}{|x-z|^{3+\gamma p}}
\nonumber
\\
&&\hspace*{55pt} {}\times\mathbb{E} \biggl(\biggl| \int_0^t \!\int
_{\mathbb{R}^3} (\pi _n\circ\tau_n ) \bigl(
\bigl[G(t-\cdot,x-\ast)- G(t-\cdot,z-\ast) \bigr]\nonumber \\
&&\hspace*{71pt}\qquad\qquad{}\times B \bigl(X_n^{-,m-1}(
\cdot,\ast) \bigr) (s,y)
M(\mathrm{d}s,\mathrm{d}y)\biggr|^p1_{L_n(t)} \biggr)\qquad
\nonumber
\\[-8pt]
\label{s2.22.22}
\\[-8pt]
\nonumber
&&\quad\le C\int_{K(t)} \mathrm{d}x\int_{K(t)}\mathrm{d}z
\frac{1}{|x-z|^{3+\gamma
p}}
\nonumber
\\
&&\hspace*{62pt}\qquad {}\times\mathbb{E} \bigl( \bigl\llVert \bigl[G(t-\cdot,x-\ast )-G(t-\cdot,z-
\ast) \bigr]\nonumber\\
&&\hspace*{150pt}\qquad{}\times B \bigl(X_n^{-,m-1} (\cdot,\ast) \bigr) \bigr
\rrVert _{\mathcal{H}_t}^p 1_{L_n(t)} \bigr).
\nonumber
\end{eqnarray}

Consider the stochastic processes $\{B(X_n^{-,m-1}(t,x))1_{L_n(t)},
(t,x)\in[0,T]\times\mathbb{R}^3\}$. First,
we apply Lemma~\ref{lss2.1.1} to $g:=B$, $Z(t,x):=X_n^{-,m-1}(t,x)
1_{L_n(t)}$. Then, by the induction
hypothesis we see that the assumptions of Lemma~\ref{lss2.1.3} are
satisfied. This yields
$R_{n,m}^{3,2}(t) \le C$. Hence, we have proved
%
\begin{equation}
\label{s2.21}
R_{n,m}^3(t)\le C.
\end{equation}

To study the term $R_n^5(t)$, we first apply Cauchy--Schwarz inequality
to obtain
\begin{eqnarray*}
R_{n,m}^5(t)&\le & \Vert h\Vert_{\mathcal{H}_t}^p
\mathbb{E} \bigl( \bigl\llVert \bigl\llVert G(t-\cdot,\bullet-\ast )D
\bigl(X_n^{m-1}(\cdot,\ast) \bigr) \bigr\rrVert
_{\mathcal{H}_t} \bigr\rrVert _{W^{\gamma,p}(K(t))}^p1_{L_n(t)}
\bigr)
\\
&\le &  C \bigl(R_{n,m}^{5,1}(t)+R_{n,m}^{5,2}(t)
\bigr),
\end{eqnarray*}
with
\begin{eqnarray*}
R_{n,m}^{5,1}(t) &=& \mathbb{E} \bigl( \bigl\llVert \bigl\llVert
G(t-\cdot ,\bullet -\ast)D \bigl(X_n^{m-1}(\cdot,\ast)
\bigr) \bigr\rrVert _{\mathcal{H}_t} \bigr\rrVert _{L^p(K(t))}^p1_{L_n(t)}
\bigr),
\\
R_{n,m}^{5,2}(t) &=& \mathbb{E} \bigl( \bigl\llVert \bigl\llVert
G(t-\cdot ,\bullet -\ast)D \bigl(X_n^{m-1}(\cdot,\ast)
\bigr) \bigr\rrVert _{\mathcal{H}_t} \bigr\rrVert _{\gamma
,p,K(t)}^p1_{L_n(t)}
\bigr).
\end{eqnarray*}

Notice that $R_{n,m}^{5,1}(t)$ is similar to the right-hand side of
(\ref{r31}), with $B(X_n^{-,m-1})$
replaced by $D(X_n^{m-1})$. Then, by analogue arguments as for
$R_{n,m}^{3,1}(t)$, we obtain
$R_{n,m}^{5,1}(t)\le C$.


Using the triangular inequality, we see that
$R_{n,m}^{5,2}(t)$ is similar to the last term in (\ref{s2.22.22}),
with $B(X_n^{-,m-1})(t,x)$ in the
latter expression replaced by $D(X_n^{m-1})(t,x)$ in the former.
Therefore, $R_{n,m}^{5,2}(t) \le C$, and
consequently,
%
\begin{equation}
\label{s2.22} R_{n,m}^5(t)\le C.
\end{equation}

By writing explicitly the convolution operator and then using
Minkowski's inequality, we have
\[
R_{n,m}^6(t)\le C\mathbb{E} \biggl(\int
_0^t \mathrm{d}s \int_{\mathbb{R}^3} G(t-s,
\mathrm{d}y) \bigl\llVert b \bigl(X_n^{m-1}(s,\bullet-y) \bigr) \bigr
\rrVert _{W^{\gamma,p}(K(t))} 1_{L_n(s)} \biggr)^p.
\]
Next, we apply H\"older's inequality with respect to the finite measure
on $[0,T]\times\mathbb{R}^3$ given
by $\mathrm{d}s G(t-s,\mathrm{d}y)$ along with (\ref{s2.38}), (\ref{s2.40.1}) with
$g:=b$ and $Z(t,x):=X_n^{m-1}(t,x)1_{L_n}(t)$.
We obtain
\[
R_{n,m}^6(t)\le C_1+C_2 \int
_0^t \mathrm{d}s \int_{\mathbb
{R}^3}G(t-s,\mathrm{d}y)
\mathbb{E} \bigl( \bigl\llVert X_n^{m-1}(s,\bullet) \bigr
\rrVert _{W^{\gamma,p}(K(s))}^p 1_{L_n(s)} \bigr).
\]
Since $\sup_{s\in[0,T]}\int_{\mathbb{R}^3} G(s,\mathrm{d}y)<\infty$, this yields
%
\begin{equation}
\label{s2.23} R_{n,m}^6(t)\le C_1+C_2
\int_0^t \mathrm{d}s \mathbb{E} \bigl( \bigl\llVert
X_n^{m-1}(s) \bigr\rrVert _{W^{\gamma,p}(K(s))}^p
1_{L_n(s)} \bigr)\le C.
\end{equation}

We now study the contribution of $R_{n,m}^4(t)$. As we did for
$R_n^5(t)$, first we apply Cauchy--Schwarz
inequality along with (\ref{s2.10}) to obtain
%
\begin{eqnarray}
 && R_{n,m}^4(t) \nonumber\\
\label{s2.24}
 &&\quad\le   C n^{{3p}/{2}}
2^{{np}/{2}}
\\
\nonumber
&&\qquad{}\times\mathbb{E}\bigl(\bigl\llVert \bigl
\llVert G(t-\cdot,\bullet-\ast) \bigl[B \bigl(X_n^{m-1}(
\cdot, \ast) \bigr)-B \bigl(X_n^{-,m-1} (\cdot,\ast) \bigr)
\bigr]   \bigr\rrVert _{\mathcal
{H}_t}
\bigr\rrVert ^p_{W^{\gamma,p}(K(t))} 1_{L_n(t)}\bigr).
\end{eqnarray}
There are two contributions coming from the right-hand side of (\ref
{s2.24}) -- the $L^p$ norm and the
fractional norm. They will be studied separately (see the terms below
denoted by $R_{n,m}^{4,1}(t)$,
$R_{n,m}^{4,2}(t)$, resp.).

We start with the contribution of the $L^p$ norm.
From Fubini's theorem and the Lipschitz continuity of $B$, it follows that
\begin{eqnarray*}
R_{n,m}^{4,1}(t) &: =&  C n^{{3p}/{2}} 2^{{np}/{2}}
\\
&&{}\times\mathbb{E} \biggl( \biggl(\int_{K(t)} \mathrm{d}x \bigl\llVert
G(t-\cdot ,x-\ast) \bigl[B \bigl(X_n^{m-1}(\cdot,\ast)
\bigr)- B \bigl(X_n^{-,m-1}(\cdot,\ast) \bigr) \bigr] \bigr
\rrVert _{\mathcal{H}_t}^p \biggr) 1_{L_n(t)} \biggr)
\\
&\le&  C n^{{3p}/{2}} 2^{{np}/{2}}
\\
&&{}\times \int_{K(t)} \mathrm{d}x \mathbb{E} \biggl( \biggl(\int
_0^t \mathrm{d}s \bigl\llVert G(t-s,x-\ast)\\
&&\qquad\hspace*{73pt}{}\times \bigl[B
\bigl(X_n^{m-1}(s,\ast) \bigr)- B \bigl(X_n^{-,m-1}(s,
\ast) \bigr) \bigr] \bigr\rrVert _{\mathcal{H}}^2
\biggr)^{{p}/{2}}
1_{L_n(t)} \biggr)
\\
&\le &  C n^{{3p}/{2}} 2^{{np}/{2}}\sup_{(t,x)\in[0,T]\times
\mathbb{R}^3}\mathbb{E}
\bigl(\bigl\vert X_n^{m-1}(t,x)- X_n^{-,m-1}(t,x)\bigr\vert^p 1_{L_n(t)} \bigr).
\end{eqnarray*}
By Lemma~\ref{lss2.1.4},
$\sup_{m} \sup_{t\in[0,T]} R_{n,m}^{4,1}(t) \le C n^{3p} 2^{-np
[({2-\beta})/{2} ]}$.
Since $\beta\in(0,2)$, this implies
%
\begin{equation}
\label{s2.26} \sup_{n,m}\sup_{t\in[0,T]}
R_{n,m}^{4,1}(t) \le C.
\end{equation}
%

Next, we study the contribution of the fractional norm of the
right-hand side of (\ref{s2.24}):
\begin{eqnarray*}
R^{4,2}_{n,m}(t) &:=& C n^{{3p}/{2}} 2^{{np}/{2}}
\\
&&{}\times\mathbb{E} \bigl( \bigl\llVert \bigl\llVert G(t-\cdot,\bullet -\ast)
\bigl[B \bigl(X_n^{m-1}(\cdot,\ast) \bigr)-B
\bigl(X_n^{-,m-1} (\cdot,\ast) \bigr) \bigr] \bigr\rrVert
_{\mathcal{H}_t} \bigr\rrVert ^p_{\gamma,p,K(t))} 1_{L_n(t)}
\bigr).
\end{eqnarray*}
Set
\[
Z_n^m(s,y)= \bigl[B \bigl(X_n^{m-1}(s,y)
\bigr)-B(X_n^{-,m-1}(s,y) \bigr]1_{L_n(s)}.
\]
For any $0\le s\le t$, we have
\[
\sup_{s\in[0,T]} \mathbb{E} \bigl( \bigl\llVert
Z_n^m(s) \bigr\rrVert ^p_{W^{\gamma,p}(K(s)^{t-s})}
\bigr) <\infty.
\]
Indeed, this holds for $Z_n^m(s)$ replaced by $B(X_n^{m-1}(s,\cdot
))1_{L_n(s)}$ and
$B(X_n^{-,m-1}(s,\cdot)1_{L_n(s)}$, separately by the following
arguments. We rely on the induction
assumption, and for the contribution of the $L^p$ norm, we use (\ref
{s2.38}). For the contribution of
the fractional norm, we apply (\ref{s2.39}).

Thus, (\ref{s2.43}) implies
%
\begin{eqnarray}
&&\mathbb{E} \bigl( \bigl\llVert \bigl\llVert G(t-\cdot,\bullet-
\ast ) \bigl[B \bigl(X_n^{m-1}(\cdot,\ast) \bigr)-B
\bigl(X_n^{-,m-1} (\cdot,\ast) \bigr) \bigr] \bigr\rrVert
_{\mathcal{H}_t} \bigr\rrVert ^p_{\gamma,p,K(t))} 1_{L_n(t)}
\bigr)
\nonumber
\\
\label{s2.27}
&&\quad \le C\int_0^t \mathrm{d}s\mathbb{E} \bigl( \bigl
\llVert \bigl[B \bigl(X_n^{m-1}(s) \bigr)-B
\bigl(X_n^{-,m-1}(s) \bigr) \bigr] \bigr\rrVert
^p_{W^{\gamma,p}
(K(s))}1_{L_n(s)} \bigr)
\\
\nonumber
&&\quad = C\int_0^t \mathrm{d}s\mathbb{E} \bigl( \bigl\llVert
\bigl[X_n^{m-1}(s)-X_n^{-,m-1}(s) \bigr]
\bigr\rrVert ^p_{W^{\gamma,p}
(K(s))}1_{L_n(s)} \bigr),
\end{eqnarray}
where in the last equality we have used that  $B$ is affine.

By the definition (\ref{sobolevfrac}), we see that
%
\begin{eqnarray}
&&\int_0^t \mathrm{d}s \mathbb{E} \bigl(
\bigl\llVert X_n^{m-1}(s)-X_n^{-,m-1}(s)
\bigr\rrVert ^p_{W^{\gamma,p}(K(s))}1_{L_n(s)} \bigr)
\nonumber
\\
\label{s2.28}
&&\quad \le C \biggl\{\int_0^t \mathrm{d}s \mathbb{E} \bigl(
\bigl\llVert X_n^{m-1}(s)-X_n^{-,m-1}(s)
\bigr\rrVert ^p_{\gamma,p,K(s)} 1_{L_n(s)} \bigr)
\\
\nonumber
&&\hspace*{12pt}\qquad{}+\sup_{(t,x)\in[0,T]\times\mathbb{R}^3} \mathbb{E} \bigl( \bigl\llvert
X_n^{m-1}(t,x)-X_n^{-,m-1} (t,x) \bigr
\rrvert ^p \bigr) \biggr\}.
\end{eqnarray}
Using Proposition~\ref{pss2.1.7} and Lemma~\ref{lss2.1.4}, the
right-hand side of this inequality is
bounded by $\tilde g_n:=C(2^{-n\bar\eta{p}/{2}}+n^{{3p}/{2}}2^{-np({3-\beta})/{2}})$, with $\bar\eta>1$.
Notice that
\[
\sup_n n^{{3p}/{2}} 2^{{np}/{2}} \tilde
g_n \le C.
\]
With all these results, we conclude
$\sup_{n,m} \sup_{t\in[0,T]} R_{n,m}^{4,2}(t) \le C$.
Along with (\ref{s2.26}), this yields
%
\begin{equation}
\label{s2.33} \sup_{n,m}\sup_{t\in[0,T]}
R_{n,m}^4(t)\le C.
\end{equation}

Bringing together (\ref{s2.18}), (\ref{s2.20}), (\ref{s2.21}),
(\ref{s2.22}), (\ref{s2.23}),
(\ref{s2.33}), we obtain
%
\begin{equation}
\label{s2.34} \mathbb{E} \bigl( \bigl\llVert X_n^m(t)
\bigr\rrVert ^p_{W^{\gamma
,p}(K(t))}1_{L_n(t)} \bigr)\le C.
\end{equation}
By the same arguments, and using that $t_n\le t$, we also have
%
\begin{equation}
\label{s2.35} \mathbb{E} \bigl( \bigl\llVert X_n^{-,m}(t)
\bigr\rrVert ^p_{W^{\gamma
,p}(K(t))}1_{L_n(t)} \bigr)\le C.
\end{equation}

From (\ref{s2.34}), (\ref{s2.35}), we obtain (\ref{s2.19}).
\end{pf*}

\section{Auxiliary results}
\label{sa}

This section gathers some technical results that are used throughout
the paper.


\begin{lemma}
\label{lss2.1.2}
Consider the Picard iterations defined in (\ref{s2.17}), (\ref
{s2.170}), respectively.
Let $p\in[1,\infty)$. For any $n\ge1$, $m\ge0$,
%
\begin{equation}
\label{s2.40} \sup_{(t,x)\in[0,T]\times\mathbb{R}^3} E \bigl( \bigl[ \bigl\llvert
X_n^m(t,x) \bigr\rrvert ^p + \bigl\llvert
X_n^{-,m}(t,x) \bigr\rrvert ^p
\bigr]1_{L_n(t)} \bigr)\le C,
\end{equation}
where the constant $C$ does not depend on $n$, $m$.

Consequently,
%
\begin{eqnarray}
\label{s2.41}
\sup_{n\ge1}\sup_{(t,x)\in[0,T]\times\mathbb{R}^3} E \bigl( \bigl[
\bigl\llvert X_n(t,x) \bigr\rrvert ^p + \bigl\llvert
X_n^{-}(t,x) \bigr\rrvert ^p
\bigr]1_{L_n(t)} \bigr)&< &  \infty,
\\
\label{s2.36}
\sup_{n\in\mathbb{N}} \sup_{m\in\mathbb{N}}\sup
_{t\in[0,T]} \mathbb{E} \bigl( \bigl[\bigl\Vert X_n^m(t)
\bigr\Vert_{L^p(K(t))} +\bigl\Vert X_n^{-,m}(t)
\bigr\Vert_{L^p(K(t))} \bigr] 1_{L_n(t)} \bigr) &< & \infty,
\\
\label{s2.37}
\sup_{n\in\mathbb{N}} \sup_{t\in[0,T]} \mathbb{E} \bigl(
\bigl[\bigl\Vert X_n(t)\bigr\Vert_{L^p(K(t))} +\bigl\Vert X_n^{-}(t)
\bigr\Vert_{L^p(K(t))} \bigr]\ 1_{L_n(t)} \bigr) &<& \infty.
\end{eqnarray}
\end{lemma}

\begin{pf}
To establish (\ref{s2.40}), we follow the arguments of the proof of
(4.9) in \cite{Delgado--Sanz-Sole012}
with $X_n(t,x)$, $X_n^-(t,x)$ in this reference replaced by
$X_n^m(t,x)$, $X_n^{-,m}(t,x)$, respectively,
and we use induction on $m$.

For the proof of (\ref{s2.41}), we use the convergences
%
\begin{eqnarray}
 \lim_{m\to\infty}\sup_{n\ge1}\sup
_{(t,x)\in[0,T]\times\mathbb
{R}^3}\mathbb{E} \bigl( \bigl\llvert X_n^m(t,x)-X_n(t,x)
\bigr\rrvert ^p1_{L_n(t)} \bigr) &=& 0,
\nonumber
\\[-8pt]
\label{s2.410}
\\[-8pt]
\nonumber
\lim_{m\to\infty}\sup_{n\ge1}\sup
_{(t,x)\in[0,T]\times\mathbb
{R}^3}\mathbb{E} \bigl( \bigl\llvert X_n^{-,m}(t,x)-X_n^-(t,x)
\bigr\rrvert ^p1_{L_n(t)} \bigr) &=& 0,
\end{eqnarray}
along with Fatou's lemma.

Property (\ref{s2.36}) follows easily from (\ref{s2.40}), and (\ref
{s2.37}) is proved by applying
(\ref{s2.36}), (\ref{s2.410}) and Fatou's lemma.\vspace*{-2pt}
\end{pf}


\begin{lemma}
\label{lss2.1.4}
Let $p\in[1,\infty)$. For any $n\ge1$, $m\ge0$,
%
\begin{equation}
\label{s2.44}
\sup_{(t,x)\in[0,T]\times\mathbb{R}^3}\mathbb{E} \bigl(\bigl\vert
X_n^{m}(t,x)-X_n^{-,m}(t,x) \bigr)
\bigr\vert^p 1_{L_n(t)} )\le C n^{{3p}/{2}} 2^{-np({3-\beta})/{2}},
\end{equation}
where the constant $C$ does not depend neither on $n$ nor on $m$. Consequently,
%
\begin{equation}
\label{s2.45}
\sup_{(t,x)\in[0,T]\times\mathbb{R}^3}\mathbb{E} \bigl(\bigl\vert
X_n(t,x)-X_n^{-}(t,x) \bigr)
\bigr\vert^p 1_{L_n(t)} )\le C n^{{3p}/{2}} 2^{-np({3-\beta})/{2}}.
\end{equation}
\end{lemma}

\begin{pf}
To establish (\ref{s2.44}), we follow the arguments of the proof of
(4.10) in \cite{Delgado--Sanz-Sole012}
with $X_n(t,x)$, $X_n^-(t,x)$ in this reference replaced by
$X_n^m(t,x)$, $X_n^{-,m}(t,x)$, respectively, and
we use induction on $m$. Then we obtain (\ref{s2.45}) by applying
(\ref{s2.410}) and Fatou's lemma.
\end{pf}

In the next lemma, we establish some results that have been shown in
the proof of \cite{Dalang--Sanz-Sole09}, Theorem~4.6, in a particular
context.

\begin{lemma}
\label{lss2.1.1}
Let $p\in[1,\infty)$, $\gamma\in(0,1)$. Consider a measurable
stochastic process $Z=\{Z(t,x), (t,x)\in[0,T]\times\mathbb{R}^3\}$
such that
\[
\sup_{(t,x)\in[0,T]\times\mathbb{R}^3} \mathbb{E} \bigl( \bigl\llvert Z(t,x) \bigr\rrvert
^p \bigr) < \infty.
\]
Let $g$ be a real-valued Lipschitz continuous function, and $K\subset
\mathbb{R}^3$ a compact set. The following properties are satisfied:
\begin{longlist}[(iii)]
\item[(i)]\vspace*{-6pt}
%
\begin{equation}
\label{s2.38}
\sup_{t\in[0,T]} \mathbb{E} \bigl(\bigl\Vert g \bigl(Z(t)
\bigr)\bigr\Vert ^p_{L^p(K)} \bigr)\le C \Bigl(1+\sup
_{t\in[0,T]}\mathbb{E} \bigl(\bigl\Vert Z(t)\bigr\Vert ^p_{L^p(K)}
\bigr) \Bigr).
\end{equation}
\item[(ii)] For any $0\le s\le t\le T$,
%
\begin{equation}
\label{s2.39}
\mathbb{E} \bigl(\bigl\Vert g \bigl(Z(s) \bigr)\bigr\Vert_{\gamma,p,K(t)^{t-s}}^p
\bigr) \le C \mathbb{E} \bigl(\bigl\Vert Z(s)\bigr\Vert^p_{\gamma,p,K(s)}
\bigr).
\end{equation}
\item[(iii)] For any $0\le s\le t\le T$, $|y|\le t-s$,
%
\begin{equation}
\label{s2.40.1} \mathbb{E} \bigl(\bigl\Vert g \bigl(Z(s, \bullet-y) \bigr)
\bigr\Vert^p_{\gamma
,p,K(t)} \bigr) \le C \mathbb{E} \bigl(\bigl\Vert Z(s)
\bigr\Vert^p_{\gamma,p,K(s)} \bigr).
\end{equation}
\end{longlist}
\end{lemma}

\begin{pf}
The assertion (i) follows by using the linear growth of the
function $g$.

Fix $ 0\le s\le t\le T$. Since $K(t)^{t-s}\subset K(s)$, we have
%
\begin{eqnarray}
\mathbb{E} \bigl(\bigl\Vert g \bigl(Z(s) \bigr)\bigr\Vert_{\gamma,p,K(t)^{t-s}}^p
\bigr) & =& \mathbb{E} \biggl(\int_{K(t)^{t-s}} \mathrm{d}x \int
_{K(t)^{t-s}} \mathrm{d}y \frac
{|g(Z(s,x))-g(Z(s,y))|^p}{|x-y|^{3+\gamma p}} \biggr)
\nonumber
\\
\label{s2.4.40}
& \le &  C \mathbb{E} \biggl(\int_{K(t)^{t-s}} \mathrm{d}x \int
_{K(t)^{t-s}} \mathrm{d}y \frac
{|Z(s,x)-Z(s,y)|^p}{|x-y|^{3+\gamma p}} \biggr)
\\
\nonumber
& =&  C \mathbb{E} \bigl(\bigl\Vert Z(s)\bigr\Vert_{\gamma,p,K(t)^{t-s}}^p
\bigr)
\\
& \le &  C \mathbb{E} \bigl(\bigl\Vert Z(s)\bigr\Vert_{\gamma,p,K(s)}^p \bigr),
\nonumber
\end{eqnarray}
which proves (ii).

Let $ 0\le s\le t\le T$ and $|y|\le t-s $. By definition,
\[
\mathbb{E} \bigl(\bigl\Vert g \bigl(Z(s,\bullet-y) \bigr)\bigr\Vert_{\gamma,p,K(t)}^p
\bigr) = \mathbb{E} \biggl(\int_{K(t)} \mathrm{d}x \int
_{K(t)} \mathrm{d}z \frac{\llvert
g(Z(s,x-y))-g(Z(s,z-y))\rrvert ^p}{|x-z|^{3+\gamma p}} \biggr).
\]
Consider the change of variables $\bar x\mapsto x-y$ and $\bar z\mapsto
z-y$. Since $x,z\in K(t)$ and $|y|\le t-s $, then $\bar x, \bar z\in
K(s)$. Thus,
\begin{eqnarray*}
\mathbb{E} \bigl(\bigl\Vert g \bigl(Z(s,\bullet-y) \bigr)\bigr\Vert_{\gamma,p,K(t)}^p
\bigr) &\le & \mathbb{E} \biggl(\int_{K(s)} \mathrm{d}\bar x \int
_{K(s)} \mathrm{d}\bar z \frac
{\llvert  g(Z(s,\bar x))-g(Z(s,\bar z))\rrvert ^p}{|\bar x-
\bar z|^{3+\gamma p}} \biggr)
\\
& =&  \mathbb{E} \bigl(\bigl\Vert g \bigl(Z(s) \bigr)\bigr\Vert_{\gamma,p,K(s)}^p
\bigr)
\\
&\le &  C \mathbb{E} \bigl(\bigl\Vert Z(s)\bigr\Vert^p_{\gamma,p,K(s)} \bigr).
\end{eqnarray*}
\upqed\end{pf}

\begin{lemma}
\label{aux1}
Let $\eta=\inf(\frac{4-\beta}{2}, 3-2\gamma-\frac{6}{p}-\beta)$,
with $\gamma\in(0,\frac{2-\beta}{2}-\frac{3}{p})$.
Let $p>\frac{2(4-\beta)}{2-\beta}$. Then
%
\begin{equation}
\label{l1} p>\frac{2\eta}{\eta-1},
\end{equation}
equivalently $\eta_1:=\frac{\eta-1}{2}-\frac{\eta}{p}>0$.
\end{lemma}
\begin{pf}
Consider first the case $\eta=\frac{4-\beta}{2}$. Then
$\frac{2\eta}{\eta-1}=\frac{2(4-\beta)}{2-\beta}$, and the
conclusion is obvious.
Next, we suppose that $\eta=3-2\gamma-\frac{6}{p}-\beta$. Then
$\frac{2\eta}{\eta-1}=\frac{6-4\gamma-{12}/{p}-2\beta
}{2-2\gamma-{6}/{p}-\beta}$.
Fix $\gamma$. The function $f\dvtx [0,\infty)\rightarrow\mathbb{R}$
defined by
\[
f(x)= \frac{6-4\gamma-2x-2\beta}{2-2\gamma-x-\beta},
\]
is increasing and $f(x)\le\lim_{x\to\infty}f(x)=2$.
Choose $p>2$. Then
\[
p>2=\sup_{x\ge0}f(x)>\frac{2\eta}{\eta-1}.
\]
Notice that $\frac{2(4-\beta)}{2-\beta}>2$. Hence, (\ref{l1}) holds.
\end{pf}
%

\section*{Acknowledgements}

The authors would like to thank an anonymous referee for its
careful reading and useful comments that helped to improve the
first version of the paper.

The first and second author where supported by the grant MICINN-FEDER
MTM 2009-07203 from the \textit{Direcci\'on General de
Investigaci\'on, Ministerio de Educaci\'on y Ciencia}, Spain.

The second author was supported by the grant MTM 2012-31192 from the
\textit{Direcci\'on General de
Investigaci\'on, Ministerio de Econom\'{\i}a y Competitividad}, Spain.






\printhistory

\begin{thebibliography}{17}


\bibitem{aida-kusuoka-stroock93}
\begin{bincollection}[mr]
\bauthor{\bsnm{Aida},~\bfnm{S.}\binits{S.}},
\bauthor{\bsnm{Kusuoka},~\bfnm{S.}\binits{S.}} \AND
\bauthor{\bsnm{Stroock},~\bfnm{D.}\binits{D.}}
(\byear{1993}).
\btitle{On the support of {W}iener functionals}.
In \bbooktitle{Asymptotic Problems in Probability Theory: {W}iener Functionals and Asymptotics ({S}anda/{K}yoto, 1990)}.
\bseries{Pitman Res. Notes Math. Ser.}
\bvolume{284}
\bpages{3--34}.
\blocation{Harlow}:
\bpublisher{Longman Sci. Tech.}
\bid{mr={1354161}}
\end{bincollection}
%

\bptok{imsref}%
\endbibitem

\bibitem{Bally--Millet-Sanz-Sole95}
\begin{barticle}[mr]
\bauthor{\bsnm{Bally},~\bfnm{Vlad}\binits{V.}},
\bauthor{\bsnm{Millet},~\bfnm{Annie}\binits{A.}} \AND
\bauthor{\bsnm{Sanz-Sol{\'e}},~\bfnm{Marta}\binits{M.}}
(\byear{1995}).
\btitle{Approximation and support theorem in H\"older norm for parabolic stochastic partial differential equations}.
\bjournal{Ann. Probab.}
\bvolume{23}
\bpages{178--222}.
\bid{issn={0091-1798}, mr={1330767}}
\end{barticle}
%

\bptok{imsref}%
\endbibitem

\bibitem{Dalang99}
\begin{barticle}[mr]
\bauthor{\bsnm{Dalang},~\bfnm{Robert~C.}\binits{R.C.}}
(\byear{1999}).
\btitle{Extending the martingale measure stochastic integral with applications to spatially homogeneous s.p.d.e.'s}.
\bjournal{Electron. J. Probab.}
\bvolume{4}
\bpages{29 pp. (electronic)}.
\bid{doi={10.1214/EJP.v4-43}, issn={1083-6489}, mr={1684157}}
\bptnote{check pages}%
\end{barticle}
%

\bptok{imsref}%
\endbibitem

\bibitem{Dalang-Quer011}
\begin{barticle}[mr]
\bauthor{\bsnm{Dalang},~\bfnm{Robert~C.}\binits{R.C.}} \AND
\bauthor{\bsnm{Quer-Sardanyons},~\bfnm{Llu{\'{\i}}s}\binits{L.}}
(\byear{2011}).
\btitle{Stochastic integrals for spde's: A comparison}.
\bjournal{Expo. Math.}
\bvolume{29}
\bpages{67--109}.
\bid{doi={10.1016/j.exmath.2010.09.005}, issn={0723-0869}, mr={2785545}}
\end{barticle}
%

\bptok{imsref}%
\endbibitem

\bibitem{Dalang--Sanz-Sole09}
\begin{barticle}[mr]
\bauthor{\bsnm{Dalang},~\bfnm{Robert~C.}\binits{R.C.}} \AND
\bauthor{\bsnm{Sanz-Sol{\'e}},~\bfnm{Marta}\binits{M.}}
(\byear{2009}).
\btitle{H\"older--{S}obolev regularity of the solution to the stochastic wave equation in dimension three}.
\bjournal{Mem. Amer. Math. Soc.}
\bvolume{199}
\bpages{vi+70}.
\bid{doi={10.1090/memo/0931}, issn={0065-9266}, mr={2512755}}
\end{barticle}
%

\bptok{imsref}%
\endbibitem

\bibitem{Dalang--Sanz-Sole13}
\begin{barticle}[mr]
\bauthor{\bsnm{Dalang},~\bfnm{Robert~C.}\binits{R.C.}} \AND
\bauthor{\bsnm{Sanz-Sol{\'e}},~\bfnm{Marta}\binits{M.}}
(\byear{2015}).
\btitle{Hitting probabilities for systems of stochastic waves}.
\bjournal{Mem. Amer. Math. Soc.}
\bvolume{237}
\bpages{v+75}.
\end{barticle}
%

\bptok{imsref}%
\endbibitem

\bibitem{Delgado--Sanz-Sole012}
\begin{barticle}[mr]
\bauthor{\bsnm{Delgado-Vences},~\bfnm{Francisco~J.}\binits{F.J.}} \AND
\bauthor{\bsnm{Sanz-Sol{\'e}},~\bfnm{Marta}\binits{M.}}
(\byear{2014}).
\btitle{Approximation of a stochastic wave equation in dimension three, with application to a support theorem in H\"older norm}.
\bjournal{Bernoulli}
\bvolume{20}
\bpages{2169--2216}.
\bid{doi={10.3150/13-BEJ554}, issn={1350-7265}, mr={3263102}}
\end{barticle}
%

\bptok{imsref}%
\endbibitem

\bibitem{Folland76}
\begin{bbook}[mr]
\bauthor{\bsnm{Folland},~\bfnm{Gerald~B.}\binits{G.B.}}
(\byear{1976}).
\btitle{Introduction to Partial Differential Equations}.
\blocation{Princeton, NJ}:
\bpublisher{Princeton Univ. Press}.
\bid{mr={0599578}}
\end{bbook}
%

\bptok{imsref}%
\endbibitem

\bibitem{Gyongy--Nualart--Sanz-Sole95}
\begin{barticle}[mr]
\bauthor{\bsnm{Gy{\"o}ngy},~\bfnm{Istv{\'a}n}\binits{I.}},
\bauthor{\bsnm{Nualart},~\bfnm{David}\binits{D.}} \AND
\bauthor{\bsnm{Sanz-Sol{\'e}},~\bfnm{Marta}\binits{M.}}
(\byear{1995}).
\btitle{Approximation and support theorems in modulus spaces}.
\bjournal{Probab. Theory Related Fields}
\bvolume{101}
\bpages{495--509}.
\bid{doi={10.1007/BF01202782}, issn={0178-8051}, mr={1327223}}
\end{barticle}
%

\bptok{imsref}%
\endbibitem

\bibitem{kx}
\begin{barticle}[mr]
\bauthor{\bsnm{Khoshnevisan},~\bfnm{Davar}\binits{D.}} \AND
\bauthor{\bsnm{Xiao},~\bfnm{Yimin}\binits{Y.}}
(\byear{2009}).
\btitle{Harmonic analysis of additive L\'evy processes}.
\bjournal{Probab. Theory Related Fields}
\bvolume{145}
\bpages{459--515}.
\bid{doi={10.1007/s00440-008-0175-5}, issn={0178-8051}, mr={2529437}}
\end{barticle}
%

\bptok{imsref}%
\endbibitem

\bibitem{Mattila95}
\begin{bbook}[mr]
\bauthor{\bsnm{Mattila},~\bfnm{Pertti}\binits{P.}}
(\byear{1995}).
\btitle{Geometry of Sets and Measures in {E}uclidean Spaces: Fractals and Rectifiability}.
\bseries{Cambridge Studies in Advanced Mathematics}
\bvolume{44}.
\blocation{Cambridge}:
\bpublisher{Cambridge Univ. Press}.
\bid{doi={10.1017/CBO9780511623813}, mr={1333890}}
\end{bbook}
%

\bptok{imsref}%
\endbibitem

\bibitem{Millet--Sanz-Sole94a}
\begin{barticle}[mr]
\bauthor{\bsnm{Millet},~\bfnm{Annie}\binits{A.}} \AND
\bauthor{\bsnm{Sanz-Sol{\'e}},~\bfnm{Marta}\binits{M.}}
(\byear{1994}).
\btitle{The support of the solution to a hyperbolic {SPDE}}.
\bjournal{Probab. Theory Related Fields}
\bvolume{98}
\bpages{361--387}.
\bid{doi={10.1007/BF01192259}, issn={0178-8051}, mr={1262971}}
\end{barticle}
%

\bptok{imsref}%
\endbibitem

\bibitem{Millet--Sanz-Sole94b}
\begin{bincollection}[mr]
\bauthor{\bsnm{Millet},~\bfnm{Annie}\binits{A.}} \AND
\bauthor{\bsnm{Sanz-Sol{\'e}},~\bfnm{Marta}\binits{M.}}
(\byear{1994}).
\btitle{A simple proof of the support theorem for diffusion processes}.
In \bbooktitle{S\'eminaire de {P}robabilit\'es, {XXVIII}}.
\bseries{Lecture Notes in Math.}
\bvolume{1583}
\bpages{36--48}.
\blocation{Berlin}:
\bpublisher{Springer}.
\bid{doi={10.1007/BFb0073832}, mr={1329099}}
\end{bincollection}
%

\bptok{imsref}%
\endbibitem

\bibitem{Millet--Sanz-Sole00}
\begin{barticle}[mr]
\bauthor{\bsnm{Millet},~\bfnm{Annie}\binits{A.}} \AND
\bauthor{\bsnm{Sanz-Sol{\'e}},~\bfnm{Marta}\binits{M.}}
(\byear{2000}).
\btitle{Approximation and support theorem for a wave equation in two space dimensions}.
\bjournal{Bernoulli}
\bvolume{6}
\bpages{887--915}.
\bid{doi={10.2307/3318761}, issn={1350-7265}, mr={1791907}}
\end{barticle}
%

\bptok{imsref}%
\endbibitem

\bibitem{V--S-S09}
\begin{barticle}[mr]
\bauthor{\bsnm{Sanz-Sol{\'e}},~\bfnm{Marta}\binits{M.}} \AND
\bauthor{\bsnm{Vuillermot},~\bfnm{Pierre~A.}\binits{P.A.}}
(\byear{2009}).
\btitle{Mild solutions for a class of fractional {SPDE}s and their sample paths}.
\bjournal{J. Evol. Equ.}
\bvolume{9}
\bpages{235--265}.
\bid{doi={10.1007/s00028-009-0014-x}, issn={1424-3199}, mr={2511552}}
\end{barticle}
%

\bptok{imsref}%
\endbibitem

\bibitem{Shimakura92}
\begin{bbook}[mr]
\bauthor{\bsnm{Shimakura},~\bfnm{Norio}\binits{N.}}
(\byear{1992}).
\btitle{Partial Differential Operators of Elliptic Type}.
\bseries{Translations of Mathematical Monographs}
\bvolume{99}.
\blocation{Providence, RI}:
\bpublisher{Amer. Math. Soc.}
\bid{mr={1168472}}
\end{bbook}
%

\bptok{imsref}%
\endbibitem

\bibitem{stroock}
\begin{binproceedings}[mr]
\bauthor{\bsnm{Stroock},~\bfnm{Daniel~W.}\binits{D.W.}} \AND
\bauthor{\bsnm{Varadhan},~\bfnm{S.~R.~S.}\binits{S.R.S.}}
(\byear{1972}).
\btitle{On the support of diffusion processes with applications to the strong maximum principle}.
In \bbooktitle{Proceedings of the {S}ixth {B}erkeley {S}ymposium on {M}athematical {S}tatistics and {P}robability ({U}niv.
{C}alifornia, {B}erkeley, {C}alif., 1970/1971), {V}ol. III: {P}robability Theory}
\bpages{333--359}.
\blocation{Berkeley, CA}:
\bpublisher{Univ. California Press}.
\bid{mr={0400425}}
\end{binproceedings}
%

\bptok{imsref}%
\endbibitem
\end{thebibliography}
\end{document}